\documentclass[12pt]{amsart}
\usepackage{amssymb,amsmath,amsfonts,latexsym}
\usepackage{bm, enumerate}
\usepackage{hyperref, xcolor}

\newcommand{\newsection}[1]{\setcounter{equation}{0} \section{#1}}
\newcommand{\bea}{\begin{eqnarray}}
\newcommand{\eea}{\end{eqnarray}}

\newcommand{\clb}{\mathcal{B}}

\newcommand{\cld}{\mathcal{D}}

\newcommand{\clh}{\mathcal{H}}

\newcommand{\clp}{\mathcal{P}}

\newcommand{\D}{\mathbb{D}}
\newcommand{\N}{\mathbb{N}}
\newcommand{\C}{\mathbb{C}}

\def\textmatrix#1&#2\\#3&#4\\{\bigl({#1 \atop #3}\ {#2 \atop #4}\bigr)}
\def\dispmatrix#1&#2\\#3&#4\\{\left({#1 \atop #3}\ {#2 \atop #4}\right)}
\newcommand{\be}{\begin{equation}}
\newcommand{\ee}{\end{equation}}
\newcommand{\ben}{\begin{eqnarray*}}
\newcommand{\een}{\end{eqnarray*}}

\newcommand{\NI}{\noindent}

\newcommand{\bi}{\begin{itemize}}
\newcommand{\ei}{\end{itemize}}
\newcommand{\diag}{\mbox{diag}}

\newcommand\la{\langle}
\newcommand\ra{\rangle}

\newtheorem{Theorem}{\sc Theorem}[section]
\newtheorem{Lemma}[Theorem]{\sc Lemma}
\newtheorem{Proposition}[Theorem]{\sc Proposition}
\newtheorem{Corollary}[Theorem]{\sc Corollary}
\newtheorem{Definition}[Theorem]{\sc Definition}
\newtheorem{Question}{\sc Question}
\newtheorem{ass}[Theorem]{\sc Assumption}

\theoremstyle{definition}
\newtheorem{Example}[Theorem]{\sc Example}
\newtheorem{Remark}[Theorem]{\sc Remark}
\newtheorem{Remarks}[Theorem]{\sc Remarks}
\newtheorem{Note}[Theorem]{\sc Note}

\newcommand{\bt}{\begin{Theorem}}
\def\beginlem{\begin{Lemma}}
\def\beginprop{\begin{Proposition}}
\def\begincor{\begin{Corollary}}
\def\begindef{\begin{Definition}}
\def\beginexamp{\begin{Example}}
\def\beginrem{\begin{Remark}}
\def\beginq{\begin{Question}}
\def\beginass{\begin{ass}}
\def\beginnote{\begin{Note}}
\newcommand{\et}{\end{Theorem}}
\def\endlem{\end{Lemma}}
\def\endprop{\end{Proposition}}
\def\endcor{\end{Corollary}}
\def\enddef{\end{Definition}}
\def\endexamp{\end{Example}}
\def\endrem{\end{Remark}}
\def\endq{\end{Question}}
\def\endass{\end{ass}}
\def\endnote{\end{Note}}

\numberwithin{equation}{section}

\begin{document}

\title{On Quasinormality of compact perturbations of the isometries}

\author[Das]{Susmita Das}
\address{Indian Institute of Science, Department of Mathematics, Bangalore, 560012,
India}
\email{susmitadas@iisc.ac.in, susmita.das.puremath@gmail.com}


\date{\today}


\begin{abstract}
We study the compact perturbations of an isometry on a separable Hilbert space and provide a complete characterization on when they are quasinormal. Based on that, we present a complete classification for a rank-one perturbation of a unilateral shift of finite multiplicity to be quasinormal in the setting of the Hardy space. The result can also be generalized for a separable Hilbert space. As an application, we provide a complete characterization for quasinormality of a rank-one perturbation of the Hardy shift.
\end{abstract}

\maketitle

\newsection{Introduction}\label{sec: intro}

Let $\clh$ be a Hilbert space and $\clb(\clh)$ be the algebra of bounded linear operators on $\clh$. An operator $T\in\clb(\clh)$ is quasinormal if it commutes with $T^*T$. This is equivalent to saying $(T^*T-TT^*)T=0$. \textsf{Throughout, we assume that all the Hilbert spaces are separable, infinite dimensional and over the complex field $\C$.} Note that, normal operators are trivially quasinormal, and so are the isometries which include natural examples of non-normal operators. Recall, an operator $T\in\clb(\clh)$ is called an isometry if $T^*T=I$, where $I$ is the identity operator on $\clh$. The class of quasinormal operators were first introduced and studied by A. Brown in \cite{A. Brown}. As we have seen, the class contains the normal operators and isometries, the two most significant objects in operator theory. A canonical representation of quasinormal operators is given in \cite{A. Brown}, which laid foundation for all the subsequent developments in the corresponding theories.

In parallel, the operators, that are compact perturbations of an isometry or a normal operator, attracted a lot of attention in recent years. There have  been several attempts in studying their properties like hyponormality, contractivity, or to find a characterization of their invariant or hyperinvariant subspaces; \cite{A. Chattopadhyay A. Skripka}, \cite{Das-Sarkar 1}, \cite{Das-Sarkar}, \cite{Foias-Jung-Ko-Percy 1}, \cite{Foias-Jung-Ko-Percy 2}, \cite{E. Ionascu}, \cite{Jung-Lee}, \cite{Ko-Lee}, \cite{Nakamura}, \cite{Nakamura 1} are a few of such. The underlying idea is to perturb an isometry or a normal operator only slightly and observe deviations in the behaviour of the deformed operators, and then find suitable criteria to obtain (or restore) the desired (or original) properties. For instance, in \cite{Nakamura}, a complete characterization is given when a rank-one perturbation of an isometry can again be a contraction or an isometry. Further, a characterization of isometries that are perturbations with compact operators of a given isometry is provided in \cite{Serban-Turcu} in terms of certain natural factorizations. In \cite{Jung-Lee}, it is shown that for a rank-one perturbation of a normal operator, hyponormality is equivalent to normality. Recall, an operator $T\in\clb(\clh)$ is a contraction if $\|Tx\|\leq\|x\|$ for all $x\in\clh$ (\cite{NF}), and hyponormal if $T^*T\geq TT^*$ (\cite{Martin-Putinar}). A characterization for the hyponormal contractions that are finite rank perturbations of unilateral shifts of finite multiplicity, is given in \cite{Das}.

In this paper, we consider compact perturbation of an isometry and investigate its quasinormality. Recall that (\cite{Halmos}), a quasinormal operator is hyponormal but the converse is not true. Also, the unilateral weighted shift with nonzero weights on the Hardy space $H^2(\D)$ is quasinormal if and only if all the weights are equal (\cite{Halmos}). With the help of this, one can show that a rank-one perturbation of an isometry may not even be hyponormal. Infact, if for nonzero $u, v\in\clh$, $u\otimes v$ defines the rank-one operator $\big(u\otimes v\big)(f)=\la f,v\ra u,$ for all $f\in\clh$, then an easy computation shows that $T= S+z\otimes 1$ on $H^2(\D)$ is not hyponormal, where $S$ is the standard Hardy shift. However, one can show that the rank-one perturbation $S-z\otimes 1$ on $H^2(\D)$ is quasinormal. \textsf{Throughout, $S$ will denote the standard unilateral shift of multiplicity one on $H^2(\D)$}, and corresponding to an orthonormal basis $\{e_n\}_{n\geq 0}$ on $\clh$, the operator $S_k$ with $k\in \N$ will denote the unilateral shift of multiplicity $k$, defined by $S_k(e_n)=e_{n+k}$. The plan of our paper is as follows: In section 2 we study the arbitrary compact perturbations of an isometry and provide a complete characterization of their quasinormality. If $T$ is a non-isometric operator which is the sum of an isometry and a compact operator, it turns out that $(T^*T)^{1/2}$ is diagonalizable with eigenvalues converging to $1$ and each non-unital eigenvalue has finite multiplicity (see Lemma 2.1). Based on this we prove that (also see Theorem \ref{thm: main}):

\begin{Theorem}
Let $T\in\clb(\clh)$ be non-isometric and of the form $T=$ Isometry $+$ Compact. Also let $\clp$ denote the set of all non-unital eigenvalues of $(T^*T)^{1/2}$ and for $n\in\mathbb{N}$, $k_n$ denote the multiplicity of $\beta_n\in\clp$. Then $T$ is quasinormal if and only if $T$ can be written as the direct sum $T= D\oplus cU$ for $c\in\{0,1\}$, $U$ being an isometry, $D$ is a diagonal operator
$$D=\diag_{\beta_n\in\clp}\big(\beta_n e^{i t^n_j}\big)_{j=1}^{k_n}; \quad t^n_j\in\mathbb{R}, \quad j=1, \ldots, k_n,$$
and $c=0$ correspond to the zero summand on a trivial space only.
\end{Theorem}

Next, in the same section we consider completely non-unitary contractions (c.n.u. for short). Recall that (\cite{NF}), a contraction $T$ is called completely non-unitary if it has no nontrivial reducing subspace on which $T$ is unitary. By canonical decomposition (\cite{NF}), every contraction on a Hilbert space can be written as a direct sum of a unitary and a completely non-unitary contraction. Since a unitary is normal (and hence quasinormal), it is sufficient to consider c.n.u. contractions only while discussing quasinormality. The defect spaces $\cld_T, \cld_{T^*}$ of a contraction $T$ are defined as $\cld_T=\overline{\text{ran}}(I-T^*T)^{1/2}$, and $\cld_{T^*}=\overline{\text{ran}}(I-TT^*)^{1/2}$, and their dimensions are called the defect indices of $T$. We provide a complete characterization of quasinormality of c.n.u. contractions with finite defect indices \big(see Theorem \ref{cont-char}\big) in terms of certain finite rank operators which appears in our previous work (see \cite{Das}). Next we consider the case of the rank-one perturbation of the unilateral shift of arbitrary finite multiplicity on the Hardy space. If $u, v\in H^2(\D)$ are nonzero and $S_k$ denotes the unilateral shift of multiplicity $k (\in\N)$, we prove the following (see Theorem \ref{Ker-rank1}):

\begin{Theorem}
If $T=S_k+u\otimes v$ on $H^2(\D)$ is quasinormal, then $\ker(I-T^*T)$ is a finitely generated $S_k$-invariant subspace of $H^2(\D)$.
\end{Theorem}

\NI In particular, we obtain the following result:

\begin{Corollary}\label{cor shift}
Let $T=S+u\otimes v$ be on $H^2(\D)$. If $T$ is quasinormal then $\ker(I-T^*T)=\theta H^2(\D)$, where $\theta$ is either a single Blaschke factor or a product of two Blaschke factors.
\end{Corollary}

Above Theorem 1.2 plays a crucial role in obtaining explicit characterizations of quasinormality for $S_k+u \otimes v$. Note that, there are only two possibilities: the vectors $v, S_k^*u$ can either be linearly dependent or be linearly independent. In the following sections, we discuss both the cases and obtain complete sets of classifications (see Theorems \ref{rank one dependent}, and \ref{rank 2 independent}). We would like to point out that
\textsf{Theorems \ref{Ker-rank1}, \ref{rank one dependent}, \ref{rank 2 independent} can also be obtained in a general Hilbert space set up---the proof of which will follow verbatim from the existing ones, only after the obvious changes in the basis vectors.} For the sake of notational convenience, we stay on the Hardy space.

Section \ref{sec: rank1 dependent} begins with the case on $\{v, S_k^*u\}$ being linearly dependent. We prove the following theorem (see also Theorem \ref{rank one dependent}):
\begin{Theorem}
Let $T=S_k+u\otimes v$ be non-isometric on $H^2(\D)$ with $\{v, S_k^*u\}$ being linearly dependent. Then $T$ is quasinormal if and only if $v$ is an eigenvector of $S_k^*$ and
\begin{equation}\label{rk1 dependent condition}
(1+\la S_k^*u, v\ra) \la z^kv, g_i\ra+\|v\|^2\sum_{n=0}^{k-1}\la u, z^{n}\ra \la z^{n}, g_i\ra=0 \quad \forall i=1, \ldots, m (\leq k),
\end{equation}
where $\{g_1,\ldots,g_m\}$ is an orthonormal set generating $\ker(I-T^*T)$.
\end{Theorem}
Let us recall that, corresponding to $\alpha\in\D$, the function $k_{\alpha}(z)=\frac{1}{1-\bar{\alpha}z}$ in $H^2(\D)$ denotes the kernel function at $\alpha$. Then as a corollary to the above theorem, we obtain the following result (see Corollary \ref{cor shift one dependent}):
\begin{Corollary}
Let $S$ be the unilateral shift and $u, v$ be nonzero elements of $H^2(\D)$ such that $v, S^*u$ are linearly dependent and the operator $T=S+u\otimes v$ is not an isometry. Then $T$ is quasinormal if and only if there exists $\alpha\in\D$ such that $v=(1-|\alpha|^2)\overline{v(\alpha)}k_{\alpha}$ and
$$1+\overline{v(\alpha)}(1-|\alpha|^2)S^*u(\alpha)=\bar{\alpha}\overline{v(\alpha)}u(0)\text{  holds}.$$
\end{Corollary}
\NI Here we like to make the following remark.

\begin{Remark}
In \cite{Ko-Lee}, Proposition 2.5 claimed that the rank-one perturbation of the unilateral shift is not quasinormal. However, the above corollary provides a class of counter examples to this claim.
\end{Remark}

Section \ref{Quasi ind} deals with the quasinormality of $S_k+u \otimes v$ with $\{v, S_k^*u\}$ linearly independent. Theorem \ref{rank 2 independent} provides a complete characterization of their quasinormality. For $k=1$ i.e., for $T=S+u\otimes v$, two subcases arise depending on the inner function $\theta$ associated to the $\ker(I-T^*T)$ namely, $\theta$ can either be the square of a single Blaschke factor or a product of two distinct Blaschke factors \big(see corollary \ref{cor shift}\big). For $\alpha\in\D$, if $B_{\alpha}(z)$ denotes the Blaschke factor
$B_{\alpha}(z)=\frac{z-\alpha}{1-\bar{\alpha}z}$ for all $z\in\D$, then we call $T$ quasinormal of type I, in case $\theta(z)=B_{\alpha}^2(z)$, and quasinormal of type II, if there exist $\alpha,\beta\in\mathbb{D}$ with $\alpha\neq \beta$ such that $\theta(z)=B_{\alpha}(z)B_{\beta}(z)$. We provide complete characterizations for both types of quasinormality in Theorems \ref{type 1 quasinormal}, and \ref{type 2 quasinormal} which also serve as another set of counter examples to Proposition 2.5 in \cite{Ko-Lee}.
\newsection{Characterization for compact perturbations}\label{char}

In this section, we consider the compact perturbations of isometries on a Hilbert space $\clh$. We begin with the following easy lemma:

\begin{Lemma}\label{lem1}
Suppose $T\in \clb(\clh)$ has form $T=Isometry + Compact$. Then $(T^*T)^{1/2}$ is diagonalizable, with eigenvalues converging to $1$ and all non-unital eigenvalues have finite multiplicities.
\end{Lemma}

\begin{proof}
Let $T=W+K$, where $W$ is an isometry and $K$ is compact on $\clh$. Then
$$T^*T=(W^*+K^*)(W+K)=I+ W^*K+ K^*W+ K^*K.$$
Since $T^*T-I$ is compact and self-adjoint, by Spectral theorem, there exist an orthonormal basis $\{e_n\}_{n\geq1}\subseteq\clh$ and a unitary $U:\clh\rightarrow\clh$ such that
\begin{equation}\label{diagonal}
U(I-T^*T)U^*=D,
\end{equation}
where $D$ is a diagonal operator defined by $D(e_n)=\alpha_n e_n$, $\alpha_n\in\C$ for all $n$ and the sequence $\{\alpha_n\}_{n\geq 1}$ converges to zero.
Now by \eqref{diagonal} $U(T^*T)^{1/2}U^*=(I-D)^{1/2}$, and $(I-D)^{1/2}e_n=(1-\alpha_n)^{1/2}e_n$ for all $n\geq 1$. Hence $(1-\alpha_n)^{1/2}$ are eigenvalues of $(I-D)^{1/2}$, and clearly the sequence $\{(1-\alpha_n)^{1/2}\}_{n\geq 1}$ converges to $1$ as $n\rightarrow \infty$. Suppose $\alpha_k\in \{\alpha_n\}_{n\geq 1}$ be a nonzero element for some $k\geq 1$. Since $\{\alpha_n\}_{n\geq 1}$ converges to zero, $\alpha_k$ can occur only finitely many times. Hence $(1-\alpha_k)^{1/2}\neq 1$ and has finite multiplicity. Since the non-unital eigenvalues of $(T^*T)^{1/2}$ are in one-one correspondence with the nonzero $\alpha_n$'s, the proof follows.
\end{proof}
With the help of the Lemma \ref{lem1}, we now provide the following classification result:

\begin{Theorem}\label{thm: main}
Let $T\in\clb(\clh)$ be non-isometric and of the form $T=$ Isometry $+$ Compact. Next let $\clp$ denote the set of all non-unital eigenvalues of $(T^*T)^{1/2}$ and for $n\in\mathbb{N}$, $k_n$ denote the multiplicity of $\beta_n\in\clp$. Then $T$ is quasinormal if and only if $T$ can be written as the direct sum $T= D\oplus cU$ for $c\in\{0,1\}$, $U$ being an isometry, $D$ is a diagonal operator
$$D=\diag_{\beta_n\in\clp}\big(\beta_n e^{i t^n_j}\big)_{j=1}^{k_n}; \quad t^n_j\in\mathbb{R}, \quad j=1, \ldots, k_n,$$
and $c=0$ correspond to the zero summand on a trivial space only.
\end{Theorem}

\begin{proof}
Fix a non-isometric $T$ of form $T=Isometry + Compact$. Then $I-T^*T\neq 0$. Lemma \ref{lem1} yields a unitary $U:\clh\rightarrow\clh$, an orthonormal basis $\{e_n\}_{n\geq 1}\subseteq\clh$, sequence $\alpha_n\in \mathbb{C}$ convergent to $0$ and Diagonal operator $D_1$ defined by $D_1(e_n)=\alpha_ne_n$ $\forall\ n$,  such that
\begin{equation}\label{I minus D}
U(T^*T)^{1/2}U^*=(I-D_1)^{1/2}.
\end{equation}

Since $(I-T^*T)\neq 0$, $\alpha_n\neq 0$ for at least one $n$, and also any nonzero $\alpha_n$ can appear only finitely many times. Hence $\alpha_n=1$ only for finitely many $n$---which implies $\ker(I-D_1)^{1/2}$ has finite-dimension.

Let $T=V(T^*T)^{1/2}$ be the unique polar decomposition of $T$ with $V$ as the appropriate partial isometry and $\ker T=\ker V=\ker (T^*T)^{1/2}$. Then w.r.t. the same unitary $U$ above, it follows by \eqref{I minus D} $UTU^*=(UVU^*)(I-D_1)^{1/2}$. Let us denote $UTU^*=T_1$ and $UVU^*=V_1$. Then $V_1$ is a partial isometry and $T_1$ is unitarily equivalent to $T$ with the polar decomposition

\begin{equation}\label{T_1 polar decomposition}
T_1=V_1(I-D_1)^{1/2}, \quad \ker T_1=\ker V_1=\ker(I-D_1)^{1/2}.
\end{equation}

Clearly, $T$ is quasinormal iff $T_1$ is so. Let us assume that $T$ is quasinormal. Then we have

\begin{equation}\label{V D commute}
V_1(I-D_1)^{1/2}=(I-D_1)^{1/2}V_1.
\end{equation}
Two cases can arise.\\
\textsf{Case 1:} $\ker T\neq \{0\}$.
By \eqref{T_1 polar decomposition}, this is equivalent to $\ker(I-D_1)^{1/2}\neq 0$.
Without loss of generality, let $\alpha_1=\alpha_2=\ldots=\alpha_k=1$ for some $k\geq 1$ and $\alpha_n\neq 1$ for all $n\geq k+1$. Let us denote $\beta_n=(1-\alpha_n)^{1/2}$ for all $n\geq 1$, and let $E(\beta)$ denote the eigenspace of $(T_1^*T_1)^{1/2}$ corresponding to an eigenvalue $\beta$. Then by Lemma \ref{lem1}, $\clh$ can be decomposed as
\begin{equation}\label{H decompose as espaces}
\clh= E(0)\oplus \big(\oplus_{\beta_n\in\clp\setminus\{0\}}E(\beta_n)\big)\oplus E(1).
\end{equation}
Note that, $\dim{E(0)}=k$ and $\dim{E(1)}$ may be finite or infinite. For $\beta_n\in\clp\setminus\{0\}$, $\clp$ being the set of all non-unital eigenvalues of $(T^*T)^{1/2}$, let $\dim{E(\beta_n)}=k_n$, where $k_n$ is the finite multiplicity of $\beta_n$. Then without loss of generality we write $(I-D_1)^{1/2}$ as

\begin{equation}\label{diagonal decomp}
(I-D_1)^{1/2}=\diag(0,\ldots, 0)\oplus \big(\oplus_{\beta_n\in\clp\setminus\{0\}}\diag(\beta_n^1, \beta_n^2,\ldots, \beta_n^{k_n})\big)\oplus I_{E(1)},
\end{equation}
where $0$ appears in first $k$ diagonal positions, each non-unital $\beta_n\in\clp\setminus \{0\}$ occurs consecutively $k_n$ times, and $I_{E(1)}:E(1)\rightarrow E(1)$ is the identity operator.

Note by \eqref{T_1 polar decomposition} and \eqref{diagonal decomp}
\begin{equation}\label{kernel of V_1}
\ker V_1=\text{span}\{e_1,\ldots, e_k\}=E(0).
\end{equation}
We now show that the eigenspaces $E(\beta)$ reduce $V_1$ for all $\beta$. Indeed by \eqref{V D commute}, for all $n\geq 1$
\begin{equation}\label{eigenspace reduces}
\begin{split}
V_1(I-D_1)^{1/2}e_n &=(I-D_1)^{1/2}V_1 e_n \\
\implies \beta_n V_1e_n&=(I-D_1)^{1/2}V_1e_n,
\end{split}
\end{equation}
and hence $e_n\in E(\beta_n)$ implies $V_1e_n\in E(\beta_n)$. Since for eigenvalues $\beta_m, \beta_n$ with $\beta_m\neq \beta_n$ $E(\beta_m)\perp E(\beta_n)$, it follows that $E(\beta_n)$ reduces $V_1$ for all $\beta_n$.


Again, for each eigenvalue $\beta$ of $(T_1^*T_1)^{1/2}$ let us define $V_1^{E(\beta)}:E(\beta)\rightarrow E(\beta)$ by
\begin{equation}\label{small isometry decomp}
V_1^{E(\beta)}=V_1(f), \text{ if } f\in E(\beta).
\end{equation}
Then one can write via \eqref{H decompose as espaces}, \eqref{eigenspace reduces}, and \eqref{small isometry decomp}
\begin{equation}\label{isometry decomp}
V_1=V_0\oplus \big(\oplus_{\beta_n\in\clp\setminus\{0\}}V_1^{E(\beta_n)}\big)\oplus V_1^{E(1)},
\end{equation}
where $V_0:E(0)\rightarrow E(0)$ is the zero operator. Now $T_1=V_1(I-D)^{1/2}$ together with \eqref{H decompose as espaces}, \eqref{diagonal decomp}, and \eqref{isometry decomp} yield
\begin{equation}\label{T decomp}
T_1=0 \oplus \big(\oplus_{\beta_n\in\clp\setminus\{0\}}\beta_n V_1^{E(\beta_n)}\big)\oplus V_1^{E(1)},
\end{equation}
where $0:E(0)\rightarrow E(0)$ is the zero operator. For each $\beta_n\in\clp\setminus\{0\}$, $V_1^{E(\beta_n)}:E(\beta_n)\rightarrow E(\beta_n)$ is finite isometry (thereby, unitary).
So there exists a unitary $U_1^{E(\beta_n)}:E(\beta_n)\rightarrow E(\beta_n)$ with
\begin{equation}\label{small isometry diagonal}
U_1^{E(\beta_n)}V_1^{E(\beta_n)}U_1^{E(\beta_n)*}=\diag(e^{it^n_1}, e^{it^n_2},\ldots, e^{it^n_{k_n}}),\quad \text{for } t^n_j\in\mathbb{R}\text{ for all }j=1, \ldots, k_n.
\end{equation}
\begin{equation}\label{unitay for full space}
\text{Note}\quad U_1=I_{E(0)}\oplus \big(\oplus_{\beta_n\in\clp\setminus\{0\}}U_1^{E(\beta_n)}\big)\oplus I_{E(1)}
\end{equation}
defines a unitary on $\clh$ and it follows by \eqref{isometry decomp}---\eqref{unitay for full space}
\begin{equation}\label{T to diagonal}
U_1T_1U_1^*=0\oplus\big(\oplus_{\beta_n\in\clp\setminus\{0\}}\beta_n\diag(e^{it^n_1}, e^{it^n_2},\ldots, e^{it^n_{k_n}})\big)\oplus V_1^{E(1)}.
\end{equation}
\begin{equation}\label{diagal only}
D=0\oplus\big(\oplus_{\beta_n\in\clp\setminus\{0\}}\beta_n\diag(e^{it^n_1}, e^{it^n_2},\ldots, e^{it^n_{k_n}})\big).
\end{equation}
By writing the above expression, we observe the following: \textsf{1) $D$ is a diagonal operator on
$E(0)\oplus(\oplus_{\beta_n\in\clp\setminus\{0\}}E(\beta_n))$; 2) the modulus of the diagonal entries of $D$ consist of the non-unital eigenvalues of $(T^*T)^{1/2}$; 3) those eigenvalues are repeated as per their multiplicities; and 4) $V_1^{E(1)}$ is an isometry on $E(1)$.} Note when $E(1)=0$, there is no isometry part $V_1^{E(1)}$ in \eqref{T to diagonal}. In that case one can write (via \eqref{T to diagonal}) with $c$ specialized to 0 where
\begin{equation}\label{T diagonal with c}
U_1T_1U_1^*=0\oplus\big(\oplus_{\beta_n\in\clp\setminus\{0\}}\beta_n\diag(e^{it^n_1}, e^{it^n_2},\ldots, e^{it^n_{k_n}})\big)\oplus c V_1^{E(1)}.
\end{equation}
If $E(1)\neq 0$, \eqref{T diagonal with c} holds for $c=1$. As $T$ is unitarily equivalent to $T_1$, the conclusion follows.

\textsf{Case 2:}
$\ker T=\{0\}$.

In this case $E(0)=0$, and by proceeding exactly as in case 1 we can decompose $T_1$ as

\begin{equation}\label{T1 last}
T_1=\oplus_{\beta_n\in\clp}\beta_n\diag(e^{it^n_1}, e^{it^n_2},\ldots, e^{it^n_{k_n}})\oplus cV_1^{E(1)},
\end{equation}
where $c\in\{0,1\}$, and $D=\oplus_{\beta_n\in\clp}\beta_n\diag(e^{it^n_1}, e^{it^n_2},\ldots, e^{it^n_{k_n}})$ is the diagonal operator with the desired properties. Then the conclusion will follow by the same argument as in case 1.

Conversely, let $T$ be non-isometric and of the form $T=$ Isometry $+$ Compact such that it can be written as the direct sum $D\oplus cU$ for some $c\in\{0,1\}$, where $U$ is an isometry and $D$ is a diagonal operator satisfying the condition of the given statement. Since $D$ is normal and an isometry (or a zero operator) is always quasinormal, $T$ is quasinormal.
\end{proof}
The following is an easy corollary of Theorem \ref{thm: main}:
\begin{Corollary}\label{cor: reducing}
Let $T$ be a nonisometric operator on $\clh$ that can be written as $T=$ Isometry $+$ Compact. If $T$ is quasinormal then $T$ has a nontrivial, proper reducing subspace.
\end{Corollary}
\begin{proof}
If $T$ is not an isometry, $(I-T^*T)$ is a nonzero compact operator. Suppose $T$ is quasinormal. Then $T(T^*T)=(T^*T)T$ and hence $T(I-T^*T)=(I-T^*T)T$. If $\text{ran }(I-T^*T)$ is not dense, then $\overline{\text{ran}}(I-T^*T)$ is a nonzero proper reducing subspace of $T$. If $\text{ran }(I-T^*T)$ is dense in $\clh$, then following the lines of proof of Theorem \ref{thm: main}, the operator $(T^*T)^{1/2}$ has a nonzero proper finite-dimensional eigenspace that reduces $T$.
\end{proof}
We now discuss the quasinormality in the settings of a contraction operator in below subsection. As we mentioned earlier, we only consider the c.n.u. contractions, which is sufficient.

\subsection{Completely non-unitary contractions and quasinormality}\label{c.n.u. and quasinormality}

Note that, if a contraction $T$ is quasinormal, then it is also hyponormal i.e., $T^*T\geq TT^*$ (\cite{Halmos}). Again, by Douglas Lemma (\cite{Douglas}), $T^*T\geq TT^*$ is equivalent to $\cld_{T}\subseteq\cld_{T^*}$. 
However, a c.n.u. $T$ with $\cld_{T}\subseteq\cld_{T^*}$ need not be quasinormal. For example, the operator $T$ on $H^2(\D)$ defined by $T(1)=T(z)=\frac{1}{2}$ and $T(z^m)= z^{m+1}$ for all $m\geq 2$, is a c.n.u. with finite defect indices and also satisfies $\cld_{T}\subseteq\cld_{T^*}$, but not even hyponormal (see Example 3.4, \cite{Das}). Recall that \big(see Theorem 3.2, \cite{Das}\big), a c.n.u. $T$ with $\cld_{T}\subseteq\cld_{T^*}$ and $\dim\cld_{T^*}<\infty$ can be written as the sum $T=S_k+F$, where $S_k$ is the unilateral shift of multiplicity $k (= \dim(\cld_{T}\ominus\cld_{T^*}))$ and $F$ is a finite rank operator satisfying certain properties. Based on this, we now provide a characterization of their quasinormality. To begin with, we recall the Proposition 3.1 and Theorem 3.2 from \cite{Das}. We keep the same set of notions as used there:

Let us fix an orthonormal basis $\{e_m\}_{m\geq 1}$ on $\clh$. By $P_n$, we denote the orthogonal projection onto the first $n$ basis vectors $\{e_1, e_2,\ldots,e_n\}$. Corresponding to $k\geq 1$, let $S_k$ denote the unilateral shift of multiplicity $k$ on $\clh$ defined by
\begin{equation}\label{S_k}
S_k(e_m)=e_{m+k},\quad m \geq 1.
\end{equation}
With this $S_k$ and $P_n$, let us define the following family of finite rank operators $F_1, F$ and $F_r ( r\geq2)$ as follows:
We fix some scalars $x^i_{i+k}\in\C$,  $i\in \{1, 2,\cdots, n\}$ and $k\geq 1$, and define
\begin{eqnarray}\label{F F1}
 F_1(x)&=&
\begin{cases}
x_1^ie_1+x_2^ie_2+\cdots+x_{i+k}^ie_{i+k}+\cdots+x_{n+k}^ie_{n+k} & \mbox{if } i\leq n, \\
0 & \mbox{if } i\geq n+1.
\end{cases}\label{F}\\
 F&=& F_1 - S_kP_n, \label{F_1'}\\
 F_r& =& F_1^r+S_k(I-P_n)F_1^{r-1}+\cdots + S_k^{r-1}(I-P_n)F_1 \quad \forall\  r\geq 2.\label{F_r}
\end{eqnarray}

\NI Then the Proposition 3.1 in \cite{Das} states that
\begin{Proposition}\label{prop: char}
Let $n \geq 0$ and $ k \in \mathbb{N}$ be arbitrary but fixed. Let $\clh$, $\{e_m\}_{m\geq 1}$, $S_k$, $F$ and $F_r$ ($r\in \mathbb{N}$) be as above \eqref{S_k}--\eqref{F_r}. Set $T=S_k+F$, and assume that $T$ is a c.n.u. contraction with finite indices such that $\cld_T\subseteq\cld_{T^*}$, $\dim\cld_{T}=n$, and $\dim(\cld_{T^*}\ominus\cld_T)=k$. Then:
\begin{enumerate}
\item $\text{rank} (P_n-F_1^*F_1)=n$.\smallskip
\item $F_1(I-P_n)=(I-P_{n+k})F_1=0$.\smallskip
\item $\lambda\|F_1^*x\|^2-\|F_1x\|^2\leq \lambda \|P_{n+k} x\|^2-\|P_n x\|^2, \quad \forall x\in \clh$ and for some $\lambda\geq 0$.\smallskip
\item $\|F_rx\|\leq \|P_n x\|$ \text{and} $\|F_r^*x\|\leq \|P_{n+kr}x\|$  \text{hold for all} $r \in \mathbb{N}$ and $x\in \clh$,
with zero as the only common solution to the corresponding equalities.
\end{enumerate}
\end{Proposition}
Based on the above, we have the following characterization for a non-isometric c.n.u. $T$ (see also Theorem 3.2 \cite{Das}):
\begin{Theorem}\label{thm: char}
Let $T$ be a bounded linear operator on a Hilbert space $\clh$ and $n, k \in \N$. Then $T$ is a c.n.u contraction such that $\cld_T \subseteq\cld_{T^*}$ with $\dim\cld_T=n$, $\dim\cld_{T^*}<\infty$ and $\dim (\cld_{T^*}\ominus\cld_T)=k$, if and only if there exists an orthonormal basis$\{e_m\}_{m\geq 1}$ of $\clh$ with respect to which $T$ can be written as $T= S_k+F$, where $S_k$ is the unilateral shift of multiplicity $k$ and $F$ is a finite rank operator satisfying conditions $(1) - (4)$ in Proposition \ref{prop: char}.
\end{Theorem}
We are now ready to prove the following characterization result on quasinormality.
\begin{Theorem}\label{cont-char}
Let $k, n\in \mathbb{N}$ and $T$ be a c.n.u. with $\dim\cld_T=n$ and $\dim(\cld_{T^*}\ominus\cld_T)=k$. Then $T$ is quasinormal if and only if there exists an orthonormal basis $\{e_n\}_{n\geq 1}$ on $\clh$ with respect to which $T$ can be written as $T=S_k+F$ where $S_k$ and $F$ are defined as above and satisfies the following conditions:
\begin{enumerate}
\item $\text{rank }(P_n-F_1^*F_1)=n$.
\item $F_1(I-P_n)=(I-P_n)F_1=0$.
\item $\|F_1^*x\|=\|F_1x\|\quad\forall x\in\clh$.
\item $\|F_1^rx\|\leq\|P_nx\| \text{and} \|F_1^rx\|\leq \|P_{n+kr}x\|$  \text{hold for all} $r \in \mathbb{N}$ and $x\in \clh$,
with zero as the only common solution to the corresponding equalities.
\end{enumerate}
\end{Theorem}
\begin{proof}
Let $T$ be a c.n.u. with $\dim\cld_T=n$ and $\dim(\cld_{T^*}\ominus\cld_T)=k$. Assume $T$ is quasinormal. Then $\cld_T\subseteq\cld_{T^*}$ and by Theorem \ref{thm: char} there exists an orthonormal basis $\{e_m\}_{m\geq 1}$ of $\clh$ with $\text{span}\{e_1, \ldots, e_n\}=\cld_T$ such that $T$ can be written as $T=S_k+F$ where $S_k$ and $F$ are defined as in \eqref{S_k}---\eqref{F_1'} and satisfies the conditions $(1)$---$(4)$ of the Proposition \ref{prop: char}. Note that the first condition is same as that of $(1)$ in Proposition \ref{prop: char}. Since $T$ is quasinormal, $\cld_T$ reduces $T$ and hence $\text{ran } F_1\subseteq \cld_T$. Therefore $(I-P_{n+k})F_1$ is the same as $(I-P_n)F_1$ and the second property follows by the property $(2)$ in Proposition \ref{prop: char}. Note that, $T$ can now be written as
\begin{equation}\label{cnu decom}
T=F_1\oplus S_k \quad \text{ on } \cld_T\oplus (\clh\ominus\cld_T).
\end{equation}
Since $F_1$ has finite rank and $T$ is quasinormal (by assumption), $F_1$ must be normal i.e., the property $3$ in the statement follows. However, the third property can also be derived as a consequence of Theorem \ref{thm: char} and the condition $(3)$ in Proposition \ref{prop: char}. Finally by first and second condition of the statement as proved above, the equation \eqref{F_r} reduces to $F_r=F_1^r$ for all $r\in\mathbb{N}$ and as a consequence, one can deduce that $I-T^rT^{*r}=P_n-F_1^{*r}F_1^r$ (see also equation $3.19$ in \cite{Das}). Then the fourth property follows as an immediate consequence of the third property (of the statement) and the condition $(4)$ in Proposition \ref{prop: char}.

Conversely, let $T$ be a c.n.u. with $\dim\cld_T=n$ and $\dim(\cld_{T^*}\ominus\cld_T)=k$, and there exists an orthonormal basis $\{e_n\}_{n\geq 1}$ on $\clh$ with respect to which $T$ can be written as $T=S_k+F$ where $S_k$ and $F$ are defined by \eqref{S_k}--\eqref{F_r} and satisfies the conditions (1)--(4) of the statement. Note by \eqref{F F1}, one can write $T=S_k+F$ as
\begin{equation}\label{T direct sum 1}
T=F_1+S_k(I-P_n).
\end{equation}
By condition (2), it follows from \eqref{T direct sum 1} that, $T$ can be decomposed as
\begin{equation}\label{T direct sum decomp}
T=F_1\oplus S_k  \quad\text{ on } \quad \text{span}\{e_m\}_{m\leq n}\oplus \overline{\text{span}}\{e_m\}_{m\geq{n+1}}.
\end{equation}
Since $S_k$ is an isometry and $F_1$ is normal (by condition (4)), $T$ is quasinormal.
\end{proof}
\begin{Remarks}
\begin{enumerate}
\item It is easy to see that, an operator $T\in\clb(\clh)$ that can be written as $T=F\oplus S_k$, where $S_k$ is the unilateral shift of finite
multiplicity $k$ and $F$ is some finite rank operator, is quasinormal if and only if $F$ is normal. However, not every quasinormal operator $T$ can be written as a direct sum of a compact or finite rank perturbation of the unilateral shift of finite multiplicity. Any diagonal operator is such an example. The Theorem \ref{cont-char} exhibits that the class of quasinormal contractions with finite defect indices has this particular property.

\item One can show that the decomposition in Theorem \ref{cont-char} coincides with the decomposition of Theorem \ref{thm: main} when applied on a
quasinormal c.n.u. with finite indices.
\end{enumerate}
\end{Remarks}
We observed in the Corollary 3.7, \cite{Das}, a c.n.u. $T$ with $\cld_T\subseteq\cld_{T^*}$, $\dim\cld_T=1$ and $\dim\cld_{T^*}<\infty$ is always hyponormal. However $T$ need not be quasinormal. Recall by Theorem 3.6, \cite{Das}, with respect to some orthonormal basis $\{e_n\}_{n\geq 0}$ on $\clh$, $T$ can be written as $T=S_k+F$ where $S_k(e_n)=e_{n+k}$ for all $n\geq 0$, and $F$ is defined by $F(e_0)=\sum_{i=0}^{k}\alpha_ie_i$ with $\alpha_i\in\mathbb{C}$, $\sum_{i=0}^{k}|\alpha_i|^2<1$ and $F(e_n)=0$, $n\geq 1$. The matrix representation $[T]$ of $T$ with respect to $\{e_n\}_{n\geq 0}$ is given by
\begin{equation}\label{T matrix}
[T] = \begin{bmatrix}
\alpha_0 &  0 & 0 &  \cdots
\\
\vdots & \vdots & \vdots & \vdots
\\
\alpha_k & 0 & 0 & \cdots
\\
0 & 1 & 0 & \cdots
\\
0 & 0 & 1 & \cdots
\\
\vdots & \vdots & \vdots &\ddots
\end{bmatrix}.
\end{equation}
We then have the following corollary:
\begin{Corollary}\label{dim 1 quasinormal}
Let $T$ be a c.n.u. with $\cld_T\subseteq\cld_{T^*}$, $\dim\cld_T=1$, $\dim\cld_{T^*}<\infty$. Then $T$ is quasinormal if and only if $\alpha_i=0$ for all $i=1,\ldots, k$, where the scalars $\alpha_i$ are associated with the matrix representation \eqref{T matrix} of $T$.
\end{Corollary}
\begin{proof}
It is easy to verify that $\cld_T=\text{span }\{e_0\}$. Let us assume $T$ to be quasinormal. Then by (2) in Theorem \ref{cont-char}, $(I-P_1)F=0$ and hence $\alpha_i=0$ for all $i=1, \ldots, k$.

\NI Conversely, if $\alpha_i=0$ for all $i=1, \ldots, k$, then $T$ can be written as $T=\alpha_0 I_{e_0}\oplus S_k$ where $I_{e_0}:\C e_0\rightarrow \C e_{0}$ is the identity operator. It is easy to verify that the conditions (1)---(4) of Theorem \ref{cont-char} are satisfied. The proof now follows by the converse of the same theorem.
\end{proof}
We now focus on the rank-one perturbations of the unilateral shifts $S_k$ of finite multiplicities $k$. We will stick to the Hardy space set-up for the notational convenience. Except for the case of $k=1$, we will state all the results for a general Hilbert space, that can be easily proved only by making the obvious changes. We start the next subsection with the following Lemma:

\subsection{Kernel space of a rank-one perturbation of $S_k$}\label{rank-one perturbations of k-multiplicity shift}

\begin{Lemma}\label{Ker-rank1}
Let $T=S_k+u\otimes v$, where $S_k$ is the unilateral shift of multiplicity $k (\geq 1)$ and $u, v$ are nonzero elements of $H^2(\D)$. If $T$ is quasinormal, then $\ker (I-T^*T)$ is $S_k$-invariant and there exists an orthonormal set of complex functions $\{g_1,\ldots, g_m\}$, $m\leq k$ such that
$$\ker (I-T^*T)=g_1 H^2[z^k]\oplus\cdots\oplus g_m H^2[z^k],\quad where$$
$$H^2[z^k]:=\{ f(z^k): f\in H^2(\D)\}.$$
\end{Lemma}
\begin{proof}
Let $T=S_k+u\otimes v$ be quasinormal, where $S_k, u,v$ are as in the statement. If $T$ is an isometry, $\ker(I-T^*T)=H^2(\D)$ and one can write
$$H^2(\D)= H^2[z^k]\oplus zH^2[z^k]\oplus\cdots\oplus z^{k-1}H^2[z^k].$$
Assume $T$ is non-isometric. Then $(I-T^*T)$ is a nonzero finite rank operator and hence, $\overline{\text{ran}}(I-T^*T)$ is a nontivial proper closed subspace of $H^2(\D)$. Hence $\ker(I-T^*T)\neq \{0\}$. Since $T$ is quasinormal, $\overline{\text{ran}}(I-T^*T)$ reduces $T$ (follows by the proof of Corollary \ref{cor: reducing}). Hence $\ker(I-T^*T)$ is $T$-invariant. Note that
\begin{equation}\label{range of defect}
(I-T^*T)=-\big(S_k^*u\otimes v+v\otimes S_k^*u+\|u\|^2 v\otimes v\big).
\end{equation}
Hence $\overline{\text{ran}}(I-T^*T)\subseteq\text{span}\{v, S_k^*u\}$. We show $v\in \text{ran}(I-T^*T)$, proceeding in two cases:

\textsf{Case 1: } $\{v, S_k^*u\}$ is linearly dependent.

Then it follows by \eqref{range of defect}, $\text{ran}(I-T^*T)=\text{span}\{v\}$.

\textsf{Case 2:} $\{v, S_k^*u\}$ is linearly independent.

\NI Then there exist $f, g\in H^2(\D)$ such that $\la f, v\ra=0$, $\la f, S_k^*u\ra\neq 0$, and $\la g, S_k^*u\ra=0$, $\la g, v\ra\neq0$. Then by equation \eqref{range of defect},
\begin{equation}\label{defect on f}
(I-T^*T)f=-(v\otimes S_k^*u)f=-\la f, S_k^*u\ra v,
\end{equation}
\begin{equation}\label{defect on g}
(I-T^*T)g=-\big(\la g, v\ra S_k^*u+\|u\|^2\la g, v\ra v\big),
\end{equation}
and further, by \eqref{defect on f} and \eqref{defect on g} $\text{ran}(I-T^*T)=\text{span}\{v,S_k^*u\}$.

The $T$-invariance of ${\overline{\text{ran}}(I-T^*T)}^{\perp}$ implies $T$-invariance of $\text{span}\{v, S_k^*u\}^{\perp}$. Now for any $f\in\text{span}\{v, S_k^*u\}^{\perp}$,
$$Tf= (S_k+ u\otimes v)f= S_k f+\la f,v\ra u=S_kf,$$
showing that $\text{span}\{v, S_k^*u\}^{\perp}$ is $S_k$-invariant. Since $S_k$ is an isometry and $k$ is finite, by classical Wold-Kolmogorov decomposition, there exist $g_1,\ldots,g_m\in \{v, S_k^*u\}^{\perp}\ominus S_k \{v, S_k^*u\}^{\perp}$ with $m\leq k$ such that
$$\{v, S_k^*u\}^{\perp}=\ker (I-T^*T)=g_1H^2[z^k]\oplus\cdots\oplus g_mH^2[z^k],$$
where the set $\{g_1,\ldots, g_m\}$ is orthonormal and the subspace $H^2[z^k]$ is defined by
$$H^2[z^k]=\{f(z^k): f\in H^2(\D)\}.$$
\end{proof}
On a general Hilbert space, the Lemma \ref{Ker-rank1} can be reformulated as follows. The proof will be similar to that of Lemma \ref{Ker-rank1}.
\begin{Lemma}\label{ker rank1 general}
Let $T=S_k+u\otimes v$, where $S_k$ is the unilateral shift of multiplicity $k (\geq 1)$ and $u, v$ are nonzero elements on a Hilbert space $\clh$. If $T$ is quasinormal, then $\ker (I-T^*T)$ is $S_k$-invariant and there exists an orthonormal set $\{g_1,\ldots, g_m\}$ with $m\leq k$ such that
$$\ker(I-T^*T)\ominus S_k\ker(I-T^*T)=\text{span }\{g_1,\ldots, g_m\}. $$
\end{Lemma}

As an application of the Theorem \ref{Ker-rank1}, we have the following corollary on the Hardy space.
\begin{Corollary}\label{cor shift}
Let $T=S+u\otimes v$, where $u,v$ are nonzero and $S$ is the unilateral shift of multiplicity one on $H^2(\D)$. If $T$ is quasinormal then $\ker(I-T^*T)=\theta H^2(\D)$, where $\theta$ is either a single Blaschke factor or a product of two Blaschke factors.
\end{Corollary}
\begin{proof}
Let $T=S+u\otimes v$ be quasinormal. If $T$ is an isometry, $\ker(I-T^*T)=H^2(\D)$ and the corollary follows with $\theta=1$. If $T$ is not an isometry, then by Lemma \ref{Ker-rank1}, $\ker(I-T^*T)$ is a proper, closed $S$-invariant subspace of $H^2(\D)$. Hence by Beurling's theorem, there exists an inner function $\theta$ such that $\ker(I-T^*T)=\theta H^2(\D)$. Again, following the lines of the proof of lemma \ref{Ker-rank1}, $\ker(I-T^*T)=\{v, S^*u\}^{\perp}$, where $v$ and $S^*u$ are not necessarily linearly independent. So $\text{codim} \ker(I-T^*T)\leq 2$. Therefore, $\theta$ is either a single Blaschke factor, or has the form
$\theta(z)=\big(\frac{z-\alpha}{1-\bar{\alpha}z}\big)^2$, or $\theta(z)=\big(\frac{z-\beta}{1-\bar{\beta}z}\big)\big(\frac{z-\gamma}{1-\bar{\gamma}z}\big) \quad\forall z\in\D$, and for some $\alpha,\beta,\gamma\in\mathbb{C}$.
\end{proof}
We are now ready to proceed for the quasinormality of the rank one perturbation $T= S_k+u\otimes v$, $k\in\N$ on the Hardy space. Since an isometry is always quasinormal, \textsf{we will only consider the operators $S_k+u\otimes v$ that are non-isometric.} Since the vectors $v, S_k^*u$ can be either linearly dependent or linearly independent, we discuss the quasinormality in two separate theorems---Theorem \ref{rank one dependent}, and Theorem \ref{rank 2 independent} in the sections \ref{sec: rank1 dependent}, and \ref{Quasi ind} respectively.

\newsection{Quasinormality, when $\{v, S_k^*u\}$ is linearly dependent}\label{sec: rank1 dependent}

We begin this section with the quasinormal characterization of $S_k+u\otimes v$ on $H^2(\D)$ with $\{v, S_k^*u\}$ linearly dependent. As a corollary \big(see Corollary \ref{cor shift one dependent}\big), we obtain a complete classification result on quasinormality for this type of rank-one perturbations of the Hardy shift. The following theorem is the main result of this section.
\begin{Theorem}\label{rank one dependent}
Let $T=S_k+u\otimes v$ on $H^2(\D)$ be non-isometric with $\{v, S_k^*u\}$ linearly dependent. Then $T$ is quasinormal if and only if $v$ is an eigen vector of $S_k^*$ and
\begin{equation}\label{rk1 dependent condition}
(1+\la S_k^*u, v\ra) \la z^kv, g_i\ra+\|v\|^2\sum_{n=0}^{k-1}\la u, z^{n}\ra \la z^{n}, g_i\ra=0 \quad \forall i=1, \ldots, m (\leq k),
\end{equation}
where $\{g_1,\ldots,g_m\}$ is an orthonormal set associated to $\ker(I-T^*T)$ as in Lemma \ref{Ker-rank1}.
\end{Theorem}
\begin{proof}
Let $T=S_k+u\otimes v$ be on $H^2(\D)$ with $\{v, S_k^*u\}$ linearly dependent and assume $T$ to be quasinormal. Note that, $\overline{\text{ran}}(I-T^*T)=\text{span}\{v\}$ (see the proof of Lamma \ref{Ker-rank1}, case 1) and  following the proof of Corollary \ref{cor: reducing}, $\overline{\text{ran}}(I-T^*T)$ reduces $T$. Again, by Lemma \ref{Ker-rank1}, $\ker (I-T^*T)$ is $S_k$-invariant. This implies, $v$ is an eigenvector of $S_k^*$. The same lemma also says that, there exist orthonormal functions $g_,\ldots, g_m\in\ker(I-T^*T)$ with $m\leq k$ such that
\begin{equation}\label{ker v perp g}
\ker(I-T^*T)=\{v\}^{\perp}=\text{span}\{g_1,\ldots, g_m\}\oplus z^k\{v\}^{\perp}
\end{equation}
Our aim is to show that the relations in \eqref{rk1 dependent condition} hold. Since $v, S_k^*u$ are linearly dependent, there exists $r\in\mathbb{C}$ such that
\begin{equation}\label{u}
S_k^*u=rv,
\end{equation}
and further by \eqref{u}, one can write
\begin{equation}\label{u proj}
u=Pu+rz^k v,
\end{equation}
where $P$ is the orthogonal projection onto $\{1, z,\ldots, z^{k-1}\}$. Note that
\begin{equation}\label{projection of u}
Pu= u(0)+\la u, z\ra z+\cdots+\la u, z^{k-1}\ra z^{k-1},
\end{equation}
and by \eqref{u}, $r=\frac{\la S_k^*u,v\ra}{\|v\|^2}$. Hence by \eqref{u proj},
\begin{equation}\label{u calculate}
u=Pu+\frac{\la S_k^*u,v\ra}{\|v\|^2}z^kv.
\end{equation}
We now find $Pu$ and $z^kv$ in terms of $v$ and $g_i\quad (i=1,\ldots, m)$. Since $H^2(\D)=\C v\oplus  \{v\}^{\perp}$ and
$\la Pu,z^{k+l}\ra=0$ (by \eqref{projection of u}) for all $l\geq 0$, one can write by \eqref{ker v perp g}
\begin{equation}\label{Pu}
Pu=c_0v+c_1 g_1+\cdots+c_m g_m,
\end{equation}
for some $c_i\in\C,\quad i=0,1,\ldots m$. Again since $\la v,g_i\ra=0$ for all $i$ and the set $\{g_i\}_{i=1}^{m}$ is orthonormal, we have by \eqref{Pu}
\begin{equation}\label{c}
c_i=
\begin{cases}
\frac{\la Pu, v\ra}{\|v\|^2}, & \mbox{if } i=0 \\
\la Pu, g_i\ra & \mbox{otherwise}.
\end{cases}
\end{equation}
\begin{equation}\label{Pu final}
Pu=\frac{\la Pu, v\ra}{\|v\|^2}v+\sum_{i=1}^{m}\la Pu, g_i\ra g_i \quad \big(\text{by }\eqref{Pu}\text{ and }\eqref{c}\big).
\end{equation}
Similarly for $z^kv$, there exist $t_0,t_1,\ldots,t_m, s\in\C$ and $f\in\{ v\}^{\perp}$ such that
\begin{equation}\label{zkv first}
z^kv=t_0v+t_1g_1+\cdots+t_mg_m+sz^kf \quad (by\quad \eqref{ker v perp g}).
\end{equation}
Since $z^k f\in z^k \{v\}^{\perp}$, we have $\la z^k f, z^k v\ra=\la f,v\ra=0$ and also by \eqref{ker v perp g}, $\la g_i, z^kf\ra=0$ for all $i=1,\ldots, m$. These together with  \eqref{zkv first} yield
\begin{equation}\label{t}
t_i=
\begin{cases}
\frac{\la z^kv, v\ra}{\|v\|^2}, & \mbox{if } i=0 \\
\la z^kv, g_i\ra & \mbox{if } i\neq 0,
\end{cases}
\end{equation}
and  $s=0$ if $f\neq 0$. Hence it follows by \eqref{zkv first}
\begin{equation}\label{zkv final}
z^kv=\frac{\la z^kv, v\ra}{\|v\|^2}v+\sum_{i=1}^{m}\la z^kv, g_i\ra g_i.
\end{equation}
Now, $Tv=(S_k+ u\otimes v)v=z^k v+ \|v\|^2 u$ and hence
\begin{equation}\label{T on v first}
\begin{split}
Tv &=z^kv+\|v\|^2\big(Pu+\frac{\la S_k^*u,v\ra}{\|v\|^2}z^kv\big)\quad (\text{by } \eqref{u calculate}) \\
&=(1+\la S_k^*u,v\ra)z^kv+\|v\|^2Pu.
\end{split}
\end{equation}
On substituting $Pu$ and $z^kv$ (from \eqref{Pu final} and \eqref{zkv final}) in \eqref{T on v first}, a little simplification yields
\begin{equation}\label{T on v}
\begin{split}
Tv &=\Big((1+\la S_k^*u,v\ra)\frac{\la z^kv,v\ra}{\|v\|^2}+\la Pu,v\ra\Big)v+\sum_{i=1}^{m}\Big((1+\la S_k^*u,v\ra)\la z^kv,g_i\ra+\|v\|^2\la Pu,g_i\ra\Big)g_i
\end{split}
\end{equation}
Since $\text{ran}(I-T^*T) (=\C v)$ reduces $T$, it follows by \eqref{T on v}
\begin{equation}\label{necessary cond first}
(1+\la S_k^*u,v\ra)\la z^kv,g_i\ra+\|v\|^2\la Pu,g_i\ra=0 \quad\text{for all } i=1,\ldots, m.
\end{equation}
Now \eqref{projection of u} and \eqref{necessary cond first} together imply
\begin{equation}\label{necessary cond final}
(1+\la S_k^*u, v\ra) \la z^kv, g_i\ra+\|v\|^2\sum_{n=0}^{k-1}\la u, z^{n}\ra \la z^{n}, g_i\ra=0,
\end{equation}
hold for all $i=1,\ldots, m$, where  $m\leq k$.

For the converse part, let $v$ be an eigenvector of $S_k^*$ and the given relations hold. Since $v, S_k^*u$ are linearly dependent, as we noted earlier,
$\overline{\text{ran}}(I-T^*T)=\C v$. Since $S_k^*v=\lambda v$ for some $\lambda\in\C$, it follows that $\ker(I-T^*T) (=\{v\}^{\perp})$ is $S_k$-invariant and
\begin{equation}\label{eigenvector T adjoint}
T^*v=(S_k^*+v\otimes u)v=(\lambda+\la v,u\ra)v,
\end{equation}
implies that $\{v\}^{\perp}$ is $T$-invariant. Also, $T|_{\{v\}^{\perp}}=S_k|_{\{v\}^{\perp}}$. Now proceeding exactly as the first part above, the given relations (which are equivalent to\eqref{necessary cond first}) reduce the equation \eqref{T on v} to
\begin{equation}\label{sufficient T}
Tv=\Big((1+\la S_k^*u,v\ra)\frac{\la z^kv,v\ra}{\|v\|^2}+\la Pu,v\ra\Big)v.
\end{equation}
Hence on $H^2(\D) (=\C v\oplus \{v\}^{\perp})$ one can write, $T=T|_{\C v}\oplus S_k$. Since $v$ is an eigenvector of $T$ (by \eqref{sufficient T}), it follows that $T$ is quasinormal.
\end{proof}
For $k=1$ in the above theorem, if $\{v, S^*u\}$ is linearly dependent with $v=a+bz$ for some scalars $a,b\in\C$ on $H^2(\D)$, then $v$ is an eigenvector of $S^*$ if and only if $b=0$. However, $v$ is always an eigenvector of $S^{*2}$. This leads to the following corollary:
\begin{Corollary}\label{v linear parmetic kernel}
Let $T=S^2+u\otimes v$ be on $H^2(\D)$ with $\{v, S^{*2}u\}$ linearly dependent and $v=a+bz$ for some scalars $a,b$, not both zero. Suppose $T$ is not an isometry. Then $T$ is quasinormal if and only if
\begin{enumerate}
\item $\la S^{*2}u,v\ra=~-1$, \text{ and}
\item $\la u,z\ra v(0)-\la v,z\ra u(0)=0$ \text{ hold}.
\end{enumerate}
\end{Corollary}
\begin{proof}
Let $T, v, S^{*2}u$ be as in the statement, and also $T$ be non-isometric. By our previous observation, $\ker(I-T^*T)=\{v\}^{\perp}$. Clearly, $v$ is an eigenvector of $S^{*2}$. Also, it is easy to see
\begin{equation}\label{v linear polynomial}
\{v\}^{\perp}=\overline{\text{span}}\{\bar{b}-\bar{a}z, z^n: n\geq 2\}.
\end{equation}
Note that, if $\alpha=\frac{\bar{b}}{\sqrt{|a|^2+|b|^2}}$ and $\beta=\frac{\bar{a}}{\sqrt{|a|^2+|b|^2}}$, then $|\alpha|^2+|\beta|^2=1$ and $\{v\}^{\perp}$ is the parametric space $H^2_{\alpha,\beta}=\overline{\text{span }}\{\alpha+ \beta z, z^n: n\geq 2\}$ (\cite{Das 1}, \cite{Paulsen Raghupati}). It follows by \eqref{v linear polynomial}
\begin{equation}\label{v perp codim}
\{v\}^{\perp}\ominus z^2\{v\}^{\perp}=\text{span}\{\bar{b}-\bar{a}z, az^2+bz^3\},
\end{equation}
and hence the set $\{\frac{\bar{b}-\bar{a}z}{\|v\|^2}, \frac{az^2+bz^3}{\|v\|^2}\}$ forms an orthonormal basis for $\ker(I-T^*T)\ominus z^2\ker(I-T^*T)$.
By assumption, there exists a nonzero $d\in\C$ such that
\begin{equation}\label{v linear s star dependent}
S^{*2}u=d v.
\end{equation}
\begin{equation}\label{d for v}
d =\frac{\la S^{*2}u,v\ra}{\|v\|^2}\quad \big(\text{by }\eqref{v linear s star dependent}\big).
\end{equation}
And further by \eqref{v linear s star dependent}, \eqref{d for v}
\begin{equation}\label{u for v}
u(z)= u(0)+\la u,z\ra z+\frac{\la S^{*2}u,v\ra}{\|v\|^2} z^2v.
\end{equation}
In what follows, we show that the conditions (1) and (2) in the statement  are equivalent to the equations in \eqref{rk1 dependent condition} of the Theorem \ref{rank one dependent} with $k=2, m=2$, and $g_1(z)=\frac{\bar{b}-\bar{a}z}{\|v\|^2}$, $g_2(z)=\frac{az^2+bz^3}{\|v\|^2}$. Indeed, in our setting, the equations in \eqref{rk1 dependent condition} read as
\begin{equation}\label{v linear ist equ}
(1+\la S^{*2}u, v\ra) \la z^2v, \frac{\bar{b}-\bar{a}z}{\|v\|^2}\ra+\|v\|^2\sum_{n=0}^{1}\la u, z^{n}\ra \la z^{n}, \frac{\bar{b}-\bar{a}z}{\|v\|^2}\ra=0,
\end{equation}
\begin{equation}\label{v linear 2nd equ}
(1+\la S^{*2}u, v\ra) \la z^2v, \frac{az^2+bz^3}{\|v\|^2}\ra+\|v\|^2\sum_{n=0}^{1}\la u, z^{n}\ra \la z^{n}, \frac{az^2+bz^3}{\|v\|^2}\ra=0.
\end{equation}
Note that, $v(0)=a$, and $\la v,z\ra=b$. Then \eqref{v linear ist equ} simplifies to
\begin{equation}\label{v linear ist equ final}
\la u,z\ra v(0)-\la v,z\ra u(0)=0,
\end{equation}
and \eqref{v linear 2nd equ}, simplifies to
\begin{equation}\label{v linear 2nd equ final}
\la S^{*2}u,v\ra=-1.
\end{equation}
\NI The proof now follows by the Theorem \ref{rank one dependent}.
\end{proof}
The analogue of the Theorem \ref{rank one dependent} on a general Hilbert space can be stated as follows:
\begin{Theorem}\label{rank one dependent general}
Let $T\in\clb(\clh)$ be non-isometric and there exists an orthonormal basis $\{e_n\}_{n\geq 0}$ with respect to which $T$ can be written as $S_k+u\otimes v$. Suppose $v, S_k^*u$ on $\clh$ are linearly dependent. Then $T$ is quasinormal if and only if $v$ is an eigenvector of $S_k^*$ and
\begin{equation}\label{rk1 depend condition}
(1+\la S_k^*u, v\ra) \la S_kv, g_i\ra+\|v\|^2\sum_{n=0}^{k-1}\la u, e_n\ra \la e_n, g_i\ra=0 \quad \forall i=1, \ldots, m (\leq k),
\end{equation}
where $\{g_1,\ldots,g_m\}$ is an orthonormal set associated to $\ker(I-T^*T)$ as in Lemma \ref{ker rank1 general}.
\end{Theorem}

\begin{proof}
The same lines of proof with $\{e_n\}_{n\geq0}$ in place of $\{z^n\}_{n\geq0}$ in Theorem \ref{rank one dependent} will work.
\end{proof}

We now prove the following corollary corresponding to $k=1$ in Theorem \ref{rank one dependent}:

\begin{Corollary}\label{cor shift one dependent}
Let $S$ be the unilateral shift and $u, v$ be nonzero elements of $H^2(\D)$ such that $v, S^*u$ are linearly dependent and the operator $T=S+u\otimes v$ is not an isometry. Then $T$ is quasinormal if and only if there exists $\alpha\in\D$ such that $v=(1-|\alpha|^2)\overline{v(\alpha)}k_{\alpha}$ and
$$1+\overline{v(\alpha)}(1-|\alpha|^2)S^*u(\alpha)=\bar{\alpha}\overline{v(\alpha)}u(0)\text{  holds}.$$
\end{Corollary}

\begin{proof}
Let $T=S+u\otimes v$ where $S, u, v$ are as given in the statement. Let $T$ be non-isometric and quasinormal. Then by Corollary \ref{cor shift},
$\ker(I-T^*T)=\theta H^2(\D)$, where $\theta$ is a product of atmost two Blaschke factors. Since $v, S^*u$ are linearly dependent, as we noted earlier
$\overline{\text{ran}}(I-T^*T)=\C v$ \big(see case 1, Lemma \ref{Ker-rank1}\big), and hence $\ker(I-T^*T)=\{v\}^{\perp}$. This implies, there exists an $\alpha\in\D$ such that $\theta$ can be taken as $\theta(z)=\frac{z-\alpha}{1-\bar{\alpha}z}, \forall z\in\D$. Since $\ker(I-T^*T)$ has codimension $1$ and $\la \theta f, k_{\alpha}\ra=0$ for all $f\in H^2(\D)$, it follows that
\begin{equation}\label{v c alpha}
v=ck_{\alpha}, \text{ where } c (\neq 0)\in\C.
\end{equation}
Now take the inner product with $k_{\alpha}$ on both sides of \eqref{v c alpha}, to observe $c=(1-|\alpha|^2)v(\alpha)$ and
\begin{equation}\label{v final}
v=(1-|\alpha|^2)v(\alpha)k_{\alpha}.
\end{equation}
Note that, $v(\alpha)\neq 0$ and $\ker(I-T^*T)\ominus z\ker(I-T^*T)=\C\theta$. Since $T$ is quasinormal (by assumption), The conditions of Theorem \ref{rank one dependent} will hold with $k=1$ and $g_1=\theta$ i.e.,

$$(1+\la S^*u,v,\ra)\la zv,\theta\ra+\|v\|^2(\la u,1\ra \la 1,\theta\ra)=0,$$
and this further simplifies to
\begin{equation}\label{simplified shift cond}
(1+\la S^*u,v,\ra)\la v,S^*\theta\ra+\|v\|^2 u(0)\overline{\theta(0)}=0.
\end{equation}
Now, substituting $S^*\theta= (1-|\alpha|^2)\frac{1}{1-\bar{\alpha}z}$, and $v$ \big(from \eqref{v final}\big) in \eqref{simplified shift cond}, we have
\begin{equation}\label{simplified shift step2}
\Big(1+(1-|\alpha|^2)\overline{v(\alpha)}\la S^*u,k_{\alpha}\ra\Big)v(\alpha) (1-|\alpha|^2)^2 \la k_{\alpha}, \frac{1}{1-\bar{\alpha}z}\ra-\bar{\alpha}|v(\alpha)|^2(1-|\alpha|^2)u(0)=0.
\end{equation}
Since $v(\alpha)\neq 0$, the equation \eqref{simplified shift step2} further simplifies to
\begin{equation}\label{simplified shift final}
1+(1-|\alpha|^2)\overline{v(\alpha)}S^*u(\alpha)=\bar{\alpha}\overline{v(\alpha)}u(0).
\end{equation}
For the converse part, let $v=(1-|\alpha|^2)v(\alpha)k_{\alpha}$ for some $\alpha\in \D$ and the given condition is satisfied. Clearly, $v$ is an eigenvector of $S^*$. Since $v, S^*u$ are linearly dependent, it follows that $\overline{\text{ran}}(I-T^*T)=\C v$ (similar as proof of Lamma \ref{Ker-rank1}, case 1) and hence $\ker(I-T^*T)(=\{v\}^{\perp})$ is $S$-invariant with codimension $1$. Hence $\ker(I-T^*T)=\theta H^2(\D)$ where $\theta$ is a single Blaschke factor. Since $\la \theta,v\ra=0$, and $v=(1-|\alpha|^2)v(\alpha)k_{\alpha}$, it follows that $\la \theta, k_{\alpha}\ra=0$. Therefore one can take $\theta(z)=\frac{z-\alpha}{1-\bar{\alpha}z}$ for all $z\in\D$. Again,
$$\ker(I-T^*T)\ominus z\ker(I-T^*T)=\text{span}\{\theta\}.$$
Now, the given condition (same as equation \eqref{simplified shift final}) is same as the that in Theorem \ref{rank one dependent} with $k=1$ and $g_1=\theta$. Hence the proof of this part follows by the converse of Theorem \ref{rank one dependent}.
\end{proof}

\begin{Remarks}
\begin{enumerate}
\item As we mentioned earlier, Corollary \ref{cor shift one dependent} provide counterexamples for Prposition 2.5 in \cite{Ko-Lee}: Rank-one
perturbation of unilateral shift is not quasinormal.

\item Let $u,v\in H^2(\D)$ be nonzero elements with $\|u\|=1$. Suppose $T=S_k+u\otimes v$ is an isometry. Then by Proposition 1 in \cite{Nakamura},
$v=(\alpha-1)S_k^*u$ where $|\alpha|=1$ with $S_k^*u\neq 0$ and $\alpha\neq 1$. Hence \textsf{$\{v, S_k^*u\}$ must be linearly dependent if $T$ is an
isometry.} However, $v$ need not be an eigenvector of $S_k^*$. For example, $S+\frac{e^{it}-1}{3}(1+z+z^2)\otimes(1+z)$ is an isometry for $t\neq 0$ but $S^*(1+z)=1\neq\beta (1+z)$ for any $\beta\in\C$. On the other hand, \textsf{if $\{v, S_k^*u\}$ is linearly independent, then $S_k+u\otimes v$ is never an isometry} (see Proposition 1, \cite{Nakamura}). We consider quasinormality in this case in the next section.
\end{enumerate}
\end{Remarks}

\newsection{Quasinormality, when $\{v, S_k^*u\}$ is linearly independent}\label{Quasi ind}

We begin this section by setting up the following notations:
Let $P$ be the orthogonal projection onto $\{1, z,\ldots, z^{k-1}\}$. Then for $u\in H^2(\D)$, we can write
\begin{equation}\label{projection of u ind}
Pu= u(0)+\la u, z\ra z+\cdots+\la u, z^{k-1}\ra z^{k-1}.
\end{equation}
Note that
$S_k f=z^k f$ for all $f\in H^2(\D)$, and $P$ is the orthogonal projection onto $\ker S_k^*$.

\NI Let us set
\begin{eqnarray}
r_0 &=& \la S_k v, v\ra+\la u,v\ra\|v\|^2, \label{eqn1}\\
r_1 &=& \la S_k v, S_k^*u\ra+\la u,S_k^*u\ra\|v\|^2, \label{eqn2}\\
s_0 &=& (1+\la S_k^*u,v\ra)\la u,v\ra+\la Pu, v\ra, \label{eqn3}\\
s_1 &=& (1+\la S_k^*u,v\ra)\la u,S_k^*u\ra+\la Pu, S_k^*u\ra,\label{eqn4}
\end{eqnarray}
and define the operator $A:\text{span}\{v, S_k^*u\}\longrightarrow \text{span}\{v, S_k^*u\}$ by
\begin{equation}\label{define operator A}
\begin{split}
A(v)&=\big(r_0\|S_k^*u\|^2-r_1\la S_k^*u,v\ra\big)v+\big(-r_0\la v,S_k^*u\ra+r_1\|v\|^2\big)S_k^*u,\\
A(S_k^*u)&=\big(s_0\|S_k^*u\|^2-s_1\la S_k^*u,v\ra\big)v+\big(-s_0\la v,S_k^*u\ra+s_1\|v\|^2\big)S_k^*u.
\end{split}
\end{equation}
In this setting, we will consider the quasinormality of $S_k+u\otimes v$, $k\in\N$, and $v, S_k^*u$ being linearly independent. As we mentioned in the last section, such an operator is never an isometry. The following theorem characterizes their quasinormal behaviour:
\begin{Theorem}\label{rank 2 independent}
Let $T=S_k+u\otimes v$ be on $H^2(\D)$ and the vectors $v, S_k^*u$ be linearly independent. Then $T$ is quasinormal if and only if $\{v, S_k^*u\}$ is $S_k^*$-invariant, the operator $A$ in \eqref{define operator A} is normal, and for an orthonormal set $\{g_1,\ldots, g_m\}$ associated to $\ker(I-T^*T)$ (as in Lemma \ref{Ker-rank1})
\begin{eqnarray}
\la z^k v,g_i\ra+\|v\|^2 \la u,g_i\ra &=& 0, \quad \text{and} \label{rel1}\\
(1+\la S_k^*u,v\ra)\la u,g_i\ra+\sum_{n=0}^{k-1}\la u,z^n\ra\la z^n,g_i\ra &=& 0\quad  \label{rel2}
\end{eqnarray}
hold for all $i=1,\ldots, m.$
\end{Theorem}

\begin{proof}
Let $T=S_k+u\otimes v$, and $v, S_k^*u$ are linearly independent. Let us assume that $T$ is quasinormal. We will split the proof into several steps.

\NI\textsf{Step 1:} In this step we show that $\{v, S_k^*u\}$ is $S_k^*$-invariant. We begin by noting (as earlier)
$$I-T^*T=-(S_k^*u\otimes v+ v\otimes S_k^*u+ \|u\|^2 v\otimes v),$$
and also it is easy to see that \big(case 2, Lemma \ref{Ker-rank1}\big)
\begin{equation}\label{ran T independent}
\overline{\text{ran}}(I-T^*T)=\text{span}\{v, S_k^*u\}.
\end{equation}
By Lemma \ref{Ker-rank1}, $\ker(I-T^*T)$ is $S_k$-invariant and hence by \eqref{ran T independent}, $\{v, S_k^*u\}$ is $S_k^*$-invariant.

\textsf{Step 2:} In this step we establish the equations \eqref{rel1}, \eqref{rel2}.

Since $\ker(I-T^*T)$ is $S_k$-invariant, it follows by Wold-Kolmogoroff decomposition (see also Lemma \ref{Ker-rank1}), there exists an orthonormal set $\{g_1,\ldots, g_m\}\subseteq H^2(\D)$ with $m\leq k$ such that
\begin{equation}\label{ker T independent}
\{v, S_k^*u\}^{\perp}\ominus z^k\{v, S_k^*u\}^{\perp}=\text{span}\{g_1,\ldots, g_m\}.
\end{equation}
Since $T$ is quasinormal, $\overline{\text{ran }}(I-T^*T)(=\text{span }\{v, S_k^*u\})$ will reduce $T$. Note that,
\begin{equation}\label{Tv independent}
Tv= (S_k+ u\otimes v)v= z^kv+\|v\|^2u \quad \text{and}
\end{equation}
\begin{equation}\label{2nd vector first indep}
T(S_k^*u)=(S_k+ u\otimes v)S_k^*u=S_k S_k^*u+\la S_k^*u,v\ra u.
\end{equation}
As $S_kS_k^*=(I-P)$ and $P$ is the orthogonal projection onto $\text{span} \{1, \ldots, z^{k-1}\}$, by \eqref{2nd vector first indep}
\begin{equation}\label{2nd vector indep}
\begin{split}
T(S_k^*u)&=(1+\la S_k^*u, v\ra)u-Pu.
\end{split}
\end{equation}
We now express $u, z^kv, Pu$ in terms of $v, S_k^*u$ and $g_i$ for $i=1,\ldots, m$. Since
$$H^2(\D)=\text{span }\{v, S_k^*u\}\oplus \{v, S_k^*u\}^{\perp},$$ and \eqref{ker T independent} holds, it follows that there exist scalars $c_i, d_i, t_i\in\C$ with $i=0,1,\cdots, m+2$ and functions $h_0, h_1, h_2\in\{v, S_k^*u\}^{\perp}$ such that
\begin{eqnarray}
u&=& c_0 v+c_1 S_k^*u+ c_2g_1+\cdots c_{m+1}g_m+c_{m+2}z^kh_0 \label{u initial}, \\
z^kv&=&d_0 v+d_1 S_k^*u+ d_2g_1+\cdots d_{m+1}g_m+d_{m+2}z^kh_1,\label{Zkv initial} \\
Pu&=&t_0 v+t_1 S_k^*u+ t_2g_1+\cdots t_{m+1}g_m+t_{m+2}z^kh_2 \label{Pu initial}.
\end{eqnarray}
Since $v, S_k^*u$ are orthogonal to $h_0, z^kh_0$, and $g_i$ for all $i=1,\ldots, m$, it follows by \eqref{u initial}
\begin{eqnarray}
\la u,v\ra &=& c_0\|v\|^2+c_1\la S_k^*u,v\ra \label{u inner product 1},\\
\la u, S_k^*u\ra &=&c_0 \la v, S_k^*u\ra+c_1\|S_k^*u\|^2 \label{u inner product 2}\quad \text{and}
\end{eqnarray}
\begin{equation}\label{ci of u}
c_{i+1}=\begin{cases}
\la u,g_i\ra, & \mbox{if } i=1,\ldots, m, \\
0, & \mbox{if }i=m+1 \quad\text{and }  h_0\neq 0.
\end{cases}
\end{equation}
Hence by \eqref{u initial} and \eqref{u inner product 1}---\eqref{ci of u},
\begin{equation}\label{u final}
u=c_0v+c_1S_k^*u+ \sum_{i=1}^{m}\la u, g_i\ra g_i.
\end{equation}
Proceeding exactly in the same way, and using \eqref{Zkv initial}. one can check the following equalities:
\begin{eqnarray}
\la z^kv,v\ra &=& d_0\|v\|^2+d_1\la S_k^*u,v\ra \label{Zkv inner product 1},\\
\la z^kv, S_k^*u\ra &=&d_0 \la v, S_k^*u\ra+d_1\|S_k^*u\|^2 \label{Zkv inner product 2},\quad \text{and further}
\end{eqnarray}
\begin{equation}\label{di of Sku}
d_{i+1}=\begin{cases}
\la z^kv,g_i\ra & \mbox{if } i=1,\ldots, m, \\
0, & \mbox{if }i=m+1 \quad\text{and }  h_1\neq 0,
\end{cases}
\end{equation}
and equations \eqref{Zkv initial} and \eqref{Zkv inner product 1}---\eqref{di of Sku} altogether imply
\begin{equation}\label{Zkv final}
z^kv=d_0v+d_1S_k^*u+ \sum_{i=1}^{m}\la z^kv, g_i\ra g_i.
\end{equation}
Similarly, and finally by \eqref{Pu initial}
\begin{eqnarray}
\la Pu,v\ra &=& t_0\|v\|^2+t_1\la S_k^*u,v\ra \label{Pu inner product 1},\\
\la Pu, S_k^*u\ra &=&t_0 \la v, S_k^*u\ra+t_1\|S_k^*u\|^2 \label{Pu inner product 2} \quad \text{and}
\end{eqnarray}
\begin{equation}\label{ti of Pu}
t_{i+1}=\begin{cases}
\la Pu,g_i\ra, & \mbox{if } i=1,\ldots, m, \\
0, & \mbox{if }i=m+1, \text{ and }  h_2\neq 0,
\end{cases}
\end{equation}
and further by \eqref{Pu initial} and \eqref{Pu inner product 1}---\eqref{ti of Pu}
\begin{equation}\label{Pu final ind}
Pu=t_0v+t_1S_k^*u+ \sum_{i=1}^{m}\la Pu, g_i\ra g_i.
\end{equation}
We substitute values of $u$ and $z^k v$ from \eqref{u final}, \eqref{Zkv final} into \eqref{Tv independent}. A simple computation shows
\begin{equation}\label{Tv quasi step}
Tv=\big(d_0+c_0\|v\|^2\big)v+ (d_1+c_1\|v\|^2)S_k^*u+ \sum_{i=1}^{m}\Big(\la z^kv,g_i\ra+\|v\|^2\la u,g_i\ra\Big)g_i.
\end{equation}
Similarly, substituting $u, Pu$ from \eqref{u final} and \eqref{Pu final ind} to \eqref{2nd vector indep}, a simplification yields
\begin{equation}\label{T 2nd vector quasi step}
\begin{split}
T(S_k^*u) &=\Big(\big(1+\la S_k^*u,v\ra\big)c_0+t_0\Big)v+\Big(\big(1+\la S_k^*u,v\ra\big)c_1+t_1\Big)S_k^*u+ \\
& \sum_{i=1}^{m}\Big(\big(1+\la S_k^*u,v\ra\big)\la u,g_i\ra+\la Pu,g_i\ra\Big)g_i.
\end{split}
\end{equation}
Since $Tv, T(S_k^*u)\in\text{span }\{v, S_k^*u\},$ it follows by \eqref{Tv quasi step} and \eqref{T 2nd vector quasi step} that for all $i=1, \ldots, m$
\begin{equation}\label{quasi zero step1}
\la z^kv,g_i\ra+\|v\|^2\la u,g_i\ra=0, \quad \text{and}
\end{equation}
\begin{equation}\label{quasi zero step2}
\big(1+\la S_k^*u,v\ra\big)\la u,g_i\ra+\la Pu,g_i\ra=0.
\end{equation}
Note by \eqref{projection of u ind}, $Pu=\sum_{n=0}^{k-1}\la u,z^n\ra z^n$, and hence equation \eqref{quasi zero step2} further reduces to
\begin{equation}\label{quasi zero step2 final}
\big(1+\la S_k^*u,v\ra\big)\la u,g_i\ra+\sum_{n=0}^{k-1}\la u,z^n\ra \la z^n,g_i\ra=0.
\end{equation}
\textsf{Step 3:} In this step we show that the operator $A$ given in the statement is normal.

\NI Since $\text{span}\{v, S_k^*u\}$ reduces $T$, quasinormality of $T$ implies that the operator $T|_{\text{span}\{v, S_k^*u\}}$ must be normal. We show that, $A$ is a nonzero scalar multiple of $T|_{\text{span}\{v, S_k^*u\}}$.

Note that the pairs of equations \eqref{Tv quasi step},\eqref{quasi zero step1} and \eqref{T 2nd vector quasi step},\eqref{quasi zero step2} yield
\begin{eqnarray}
Tv&=&\big(d_0+c_0\|v\|^2\big)v+ (d_1+c_1\|v\|^2)S_k^*u, \text{  and} \label{Tv reduced}\\
T(S_k^*u)&=&\Big(\big(1+\la S_k^*u,v\ra\big)c_0+t_0\Big)v+\Big(\big(1+\la S_k^*u,v\ra\big)c_1+t_1\Big)S_k^*u\label{Tv reduced 2nd}
\end{eqnarray}
respectively.
We now find the scalars $c_i, d_i, t_i$ for $i=0,1$. Note that, for $v, S_k^*u$ linearly independent, Cauchy-Schwarz inequality implies
\begin{equation}\label{Cauchy-Schwarz inq}
|\la v,S_k^*u\ra|^2<\|v\|^2\|S_k^*u\|^2.
\end{equation}
With the help of \eqref{Cauchy-Schwarz inq}, it is easy to see via \eqref{u inner product 1}, \eqref{u inner product 2}
\begin{eqnarray}
c_0 &=& \frac{1}{|\la v,S_k^*u\ra|^2-\|v\|^2\|S_k^*u\|^2}\Big(\la u,v\ra\|S_k^*u\|^2-\la u,S_k^*u\ra\la S_k^*u,v\ra\Big), \label{c0}\\
c_1 &=& \frac{1}{\|\la v,S_k^*u\ra\|^2-\|v\|^2\|S_k^*u\|^2}\Big(\la u,v\ra\la v,S_k^*u\ra-\la u,S_k^*u\ra\|v\|^2\Big).\label{c1}
\end{eqnarray}
Similarly, it follows by \eqref{Zkv inner product 1}, \eqref{Zkv inner product 2}, and \eqref{Cauchy-Schwarz inq}
\begin{eqnarray}
d_0 &=& \frac{1}{|\la v,S_k^*u\ra|^2-\|v\|^2\|S_k^*u\|^2}\Big(\la z^kv,v\ra\|S_k^*u\|^2-\la z^kv,S_k^*u\ra\la S_k^*u,v\ra\Big), \label{d0}\\
d_1 &=& \frac{1}{|\la v,S_k^*u\ra|^2-\|v\|^2\|S_k^*u\|^2}\Big(\la z^kv,v\ra\la v,S_k^*u\ra-\la z^kv,S_k^*u\ra\|v\|^2\Big),\label{d1}
\end{eqnarray}
and finally by \eqref{Pu inner product 1}, \eqref{Pu inner product 2}, and \eqref{Cauchy-Schwarz inq}
\begin{eqnarray}
t_0 &=& \frac{1}{|\la v,S_k^*u\ra|^2-\|v\|^2\|S_k^*u\|^2}\Big(\la Pu,v\ra\|S_k^*u\|^2-\la Pu,S_k^*u\ra\la S_k^*u,v\ra\Big), \label{t0}\\
t_1 &=& \frac{1}{|\la v,S_k^*u\ra|^2-\|v\|^2\|S_k^*u\|^2}\Big(\la Pu,v\ra\la v,S_k^*u\ra-\la Pu,S_k^*u\ra\|v\|^2\Big).\label{t1}
\end{eqnarray}
Our aim is to simplify the coefficients of $v, S_k^*u$ in $Tv$ and $T(S_k^*u)$ in \eqref{Tv reduced}, \eqref{Tv reduced 2nd}. Substituting $c_0, d_0$ from \eqref{c0}, \eqref{d0} in $d_0+c_0\|v\|^2$, a simple computation yields
\begin{equation}\label{1st coeff of Tv}
\begin{split}
d_0+c_0\|v\|^2=& -\frac{1}{|\la v,S_k^*u\ra|^2-\|v\|^2\|S_k^*u\|^2}\Big(\big(\la z^k v, v\ra+\la u,v\ra\|v\|^2\big) \|S_k^*u\|^2- \\
&\qquad\qquad\qquad\qquad\big(\la z^k v, S_k^*u\ra+\la u,S_k^*u\ra\|v\|^2\big) \la S_k^*u, v\ra\Big)\\
&=-\frac{1}{|\la v,S_k^*u\ra|^2-\|v\|^2\|S_k^*u\|^2}\big(r_0\|S_k^*u\|^2-r_1\la S_k^*u,v\ra\big),
\end{split}
\end{equation}
where the last equality follows by using the notations $r_0, r_1$ from\eqref{eqn1}, and \eqref{eqn2} respectively.

\NI Similarly, substituting $c_1, d_1$ from \eqref{c1}, \eqref{d1} in the quantity $d_1+c_1\|v\|^2$, one will have
\begin{equation}\label{2nd coeff of Tv}
\begin{split}
d_1+c_1\|v\|^2=& \frac{1}{|\la v,S_k^*u\ra|^2-\|v\|^2\|S_k^*u\|^2}\Big(\big(\la z^k v, v\ra+\la u,v\ra\|v\|^2\big) \la v, S_k^*u\ra- \\
&\qquad\qquad\qquad\qquad\big(\la z^k v, S_k^*u\ra+\la u,S_k^*u\ra\|v\|^2\big) \|v\|^2\Big)\\
&=\frac{1}{|\la v,S_k^*u\ra|^2-\|v\|^2\|S_k^*u\|^2}\big(r_0\la v,S_k^*u\ra-r_1\|v\|^2\big) \quad (\text{ by } \eqref{eqn1}, \eqref{eqn2}).
\end{split}
\end{equation}
Next,
\begin{equation}\label{1st coeff of TSku}
\begin{split}
(1+\la S_k^*u,&\ v\ra)c_0+t_0= -\frac{1}{|\la v,S_k^*u\ra|^2-\|v\|^2\|S_k^*u\|^2}\Big(\big((1+\la S_k^*u, v\ra)\la u,v\ra+  \la Pu,v\ra\big)\|S_k^*u\|^2\\  - &\ \big((1+\la S_k^*u,  v\ra)\la  u,S_k^*u\ra+\la Pu,S_k^*u\ra\big) \la S_k^*u, v\ra\Big) \quad (\text{by }\eqref{c0}, \eqref{t0}, \eqref{Cauchy-Schwarz inq})\\
=&\ -\frac{1}{|\la v,S_k^*u\ra|^2 -\|v\|^2\|S_k^*u\|^2} \big(s_0\|S_k^*u\|^2-s_1\la S_k^*u,v\ra\big) \quad (\text{by }\eqref{eqn3}, \eqref{eqn4}),
\end{split}
\end{equation}
and finally
\begin{equation}\label{2nd coeff of TSku}
\begin{split}
(1+\la S_k^*u, v\ra)c_1+t_1=& \frac{1}{|\la v,S_k^*u\ra|^2-\|v\|^2\|S_k^*u\|^2}\Big(\big((1+\la S_k^*u, v\ra)\la u,v\ra+\la Pu,v\ra\big) \la v,S_k^*u\ra- \\
&\qquad\qquad\big((1+\la S_k^*u, v\ra)\la u,S_k^*u\ra+\la Pu,S_k^*u\ra\big) \|v\|^2\Big) \quad (\text{by }\eqref{c1}, \eqref{t1}, \eqref{Cauchy-Schwarz inq})\\
&=\frac{1}{|\la v,S_k^*u\ra|^2-\|v\|^2\|S_k^*u\|^2}\big(s_0\la v,S_k^*u\ra-s_1\|v\|^2\big) \quad (\text{by }\eqref{eqn3}, \eqref{eqn4}).
\end{split}
\end{equation}
Now Equations \eqref{1st coeff of Tv}--\eqref{2nd coeff of TSku} together with \eqref{Tv reduced},\eqref{Tv reduced 2nd}, and \eqref{define operator A} implies
$$T|_{\text{span}\{v,S_k^*u\}}=-\Big(\frac{1}{|\la v,S_k^*u\ra|^2-\|v\|^2\|S_k^*u\|^2}\Big) A.$$
For the converse part, let  $T=S_k+ u\otimes v$ be with $\{v, S_k^*u\}$ linearly independent and all the conditions in the statement hold. Then $\overline{\text{ran}}(I-T^*T)=\text{span }\{v, S_k^*u\}$ \big(see case 2, Lemma \ref{Ker-rank1}\big). Since $\{v, S_k^*u\}$ is $S_k^*$-invariant (by one of the conditions), it follows that $\{v, S_k^*u\}^{\perp} (=\ker(I-T^*T))$ is $S_k$-invariant. Also, for any $f\in \{v, S_k^*u\}^{\perp}$
$$Tf= S_k f+\la f,v\ra u= S_kf, \quad \text{and hence}$$
\begin{equation}\label{T inv ind}
T|_{\{v, S_k^*u\}^{\perp}}=S_k|_{\{v, S_k^*u\}^{\perp}}.
\end{equation}
Again, by Wold-Kolmogoroff decomposition, (see also Lemma \ref{Ker-rank1}) there exist orthonormal functions $\{g_i\}_{i=1}^m$, $m\leq k$ such that
\begin{equation}\label{same like 255}
\{v, S_k^*u\}^{\perp}\ominus z^k \{v, S_k^*u\}^{\perp}=\text{span }\{g_1,\ldots, g_m\}.
\end{equation}
As $H^2(\D)=\text{span }\{v, S_k^*u\}\oplus \{v, S_k^*u\}^{\perp}$, proceeding exactly as in the first part, \eqref{same like 255} shows
\begin{equation}\label{Tv quasi step same}
Tv=\big(d_0+c_0\|v\|^2\big)v+ (d_1+c_1\|v\|^2)S_k^*u+ \sum_{i=1}^{m}\Big(\la z^kv,g_i\ra+\|v\|^2\la u,g_i\ra\Big)g_i,\quad \text{and}
\end{equation}
\begin{equation}\label{T 2nd vector quasi step same}
\begin{split}
T(S_K^*u) &=\Big(\big(1+\la S_k^*u,v\ra\big)c_0+t_0\Big)v+\Big(\big(1+\la S_k^*u,v\ra\big)c_1+t_1\Big)S_k^*u+ \\
& \sum_{i=1}^{m}\Big(\big(1+\la S_k^*u,v\ra\big)\la u,g_i\ra+\la Pu,g_i\ra\Big)g_i
\end{split}
\end{equation}
for some scalars $c_j, d_j, t_j\in\C$ for $j=0,1$ satisfying the same set of equations \eqref{c0}---\eqref{t1}. Again, by the given conditions \eqref{rel1}, \eqref{rel2}, the equations \eqref{Tv quasi step same}, and \eqref{T 2nd vector quasi step same} reduce to
\begin{eqnarray}
Tv&=&\big(d_0+c_0\|v\|^2\big)v+ (d_1+c_1\|v\|^2)S_k^*u, \text{  and} \label{Tv reduced same}\\
T(S_k^*u)&=&\Big(\big(1+\la S_k^*u,v\ra\big)c_0+t_0\Big)v+\Big(\big(1+\la S_k^*u,v\ra\big)c_1+t_1\Big)S_k^*u.\label{Tv reduced 2nd same}
\end{eqnarray}
Clearly by \eqref{Tv reduced same}, \eqref{Tv reduced 2nd same}, $\overline{\text{ran}}(I-T^*T) (=\text{span }\{v, S_k^*u\})$ is $T$-invariant.This together with \eqref{T inv ind} implies that $T$ on $H^2(\D)\Big(=\text{span }\{v, S_k^*u\}\oplus \{v, S_k^*u\}^{\perp}$\Big) decomposes as
$$T=T|_{\text{span }\{v, S_k^*u\}}\oplus S_k.$$
Now following the same steps of simplifications as in the first part upon substituting the scalars $c_j, d_j, t_j\in\C$ for $j=0,1$ from \eqref{c0}---\eqref{t1},
one will have $$T|_{\text{span }\{v, S_k^*u\}}=-\frac{1}{|\la v,S_k^*u\ra|^2-\|v\|^2\|S_k^*u\|^2} A,$$ where $A$ is given by \eqref{define operator A}. Since $S_k$ is quasinormal, the proof now follows by the given condition: $A$ is normal.
\end{proof}
In a general Hilbert space set-up, Theorem \ref{rank 2 independent} can be rephrased as:
\begin{Theorem}\label{rank two independent general}
Let $T\in\clb(\clh)$ and there exists an orthonormal basis $\{e_n\}_{n\geq 0}$ with respect to which $T$ can be written as $S_k+u\otimes v$. Suppose $v, S_k^*u$ on $\clh$ are linearly independent. Then $T$ is quasinormal if and only if $\{v, S_k^*u\}$ is $S_k^*$-invariant, the operator $A$ in \eqref{define operator A} is normal and
\begin{eqnarray}
\la S_k v,g_i\ra+\|v\|^2 \la u,g_i\ra &=& 0,\label{rel-1}\\
(1+\la S_k^*u,v\ra)\la u,g_i\ra+\sum_{n=0}^{k-1}\la u,e_n\ra\la e_n,g_i\ra &=& 0,\label{rel-2}
\end{eqnarray}
hold for all $i=1,\ldots,m$, where $\{g_1,\ldots, g_m\}$ is an orthonormal set associated to $\ker(I-T^*T)$ as in Lemma \ref{ker rank1 general}.
\end{Theorem}
\begin{proof}
The proof is similar to that of Theorem \ref{rank 2 independent}, follows by replacing all $z^n$'s by $e_n$'s.
\end{proof}
Let us now consider the case $k=1$ in the above theorem \ref{rank 2 independent} i.e., let $T=S+ u\otimes v$ where $\{v, S^*u\}$ is linearly independent. As we noted earlier \big(see case 2, Lemma \ref{Ker-rank1}\big), $\text{ran}(I-T^*T)=\text{span}\{v, S^*u\}$ and hence $\ker(I-T^*T)=\{v, S^*u\}^{\perp}$. If $T$ is quasinormal, then it necessarily follows \big(Corollary \ref{cor shift}\big) that $\ker(I-T^*T)=\theta H^2(\D)$, where $\theta$ can be taken as either $\big(\frac{z-\alpha}{1-\bar{\alpha}z}\big)^2$ for some $\alpha\in\D$ or $\big(\frac{z-\beta}{1-\bar{\beta}z}\big)\big(\frac{z-\gamma}{1-\bar{\gamma}z}\big)$ for some distinct $\beta,\gamma\in\D$. Recall that \big(see section \ref{sec: intro}\big), $T$ is quasinormal of type I if $\ker(I-T^*T)=\big(\frac{z-\alpha}{1-\bar{\alpha}z}\big)^2$, $\alpha\in\D$, and quasinormal of type~II if $\ker(I-T^*T)=\big(\frac{z-\beta}{1-\bar{\beta}z}\big)\big(\frac{z-\gamma}{1-\bar{\gamma}z}\big)$, $\beta,\gamma\in\D$ with $\beta\neq\gamma$. With the help of the kernel functions, in the following subsections we deduce \big(from Theorem \ref{rank 2 independent}\big) more refined characterization for these two types of quasinormality.

\subsection{Type I quasinormality}\label{sec: Type 1}

\begin{Theorem}\label{type 1 quasinormal}
Let $T=S+u\otimes v$ be on $H^2(\D)$ with $\{v, S^*u\}$ linearly independent. Them $T$ is quasinormal of type I if and only if there exists $\alpha\in\D$ such that
\begin{enumerate}
\item $v, S^*u\in\text{ span } \{k_{\alpha}, \frac{z-\alpha}{1-\bar{\alpha}z}k_{\alpha}\},$
\item $v(\alpha)+\alpha B_{\alpha}^*v(\alpha)=0,$
\item $\overline{B_{\alpha}^*v(\alpha)}\la u, (\frac{z-\alpha}{1-\bar{\alpha}z})^2\ra +1=0$,
and
\item The operator $T_1:\text{ span } \{k_{\alpha}, \frac{z-\alpha}{1-\bar{\alpha}z}k_{\alpha}\}\longrightarrow \text{ span } \{k_{\alpha},
\frac{z-\alpha}{1-\bar{\alpha}z}k_{\alpha}\}$ defined by
\begin{equation*}
\begin{split}
T_1(k_{\alpha})&= \Big(\overline{v(\alpha)}u(0)+\alpha\big(1+\overline{v(\alpha)}R\big)\Big)k_{\alpha}+ \Big(1+\overline{v(\alpha)}R\Big)\frac{z-\alpha}{1-\bar{\alpha}z}k_{\alpha} \\
T_1(\frac{z-\alpha}{1-\bar{\alpha}z}k_{\alpha})&=\Big(\overline{B_{\alpha}^*v(\alpha)}u(0)+\alpha \overline{B_{\alpha}^*v(\alpha)}R\Big)k_{\alpha}+\Big(\overline{B_{\alpha}^*v(\alpha)}R\Big)\frac{z-\alpha}{1-\bar{\alpha}z}k_{\alpha}
\end{split}
\end{equation*}
is normal, where $R=(1-|\alpha|^2)S^*u(\alpha)-\bar{\alpha}u(0)$.
\end{enumerate}
\end{Theorem}

\begin{proof}
Let $T=S+ u\otimes v$ and $v, s^*u$ are linearly independent on $H^2(\D)$. As we observed earlier, $\overline{\text{ran}}(I-T^*T)=\text{ span }\{v, S^*u\}$, and hence $\ker(I-T^*T)=\{v, S^*u\}^{\perp}$. Assume that, $T$ is quasinormal of type I. Then $\text{ span}\{v, S^*u\}$ reduces $T$, and is $S^*$-invariant (see theorem \ref{rank 2 independent}). Also in addition, there exists an $\alpha\in\D$ such that $\ker(I-T^*T)=\{v, S^*u\}^{\perp}=\big(\frac{z-\alpha}{1-\bar{\alpha}z}\big)^2H^2(\D)$ (by the discussion prior to the statement). Since $k_{\alpha}, \frac{z-\alpha}{1-\bar{\alpha}z}k_{\alpha}$ are mutually orthogonal and
$$k_{\alpha}, \frac{z-\alpha}{1-\bar{\alpha}z}k_{\alpha}\perp \big(\frac{z-\alpha}{1-\bar{\alpha}z}\big)^2 z^n, \text{ for all } n\geq 0,$$
it follows that
\begin{equation}\label{v su perp kernel}
v, S^*u\in\text{ span } \{k_{\alpha}, \frac{z-\alpha}{1-\bar{\alpha}z}k_{\alpha}\},
\end{equation}
which is the condition (1) in the statement.

Note that, $\ker(I-T^*T)\ominus z\ker(I-T^*T)=\C \big(\frac{z-\alpha}{1-\bar{\alpha}z}\big)^2$, and since $T$ is quasinormal, it follows by Theorem \ref{rank 2 independent}, the relations \ref{rel1}, and \ref{rel2} hold corresponding to $g_1(z)=\big(\frac{z-\alpha}{1-\bar{\alpha}z}\big)^2$ i.e.,
\begin{eqnarray}
\la zv, \big(\frac{z-\alpha}{1-\bar{\alpha}z}\big)^2\ra+ \|v\|^2\la u, \big(\frac{z-\alpha}{1-\bar{\alpha}z}\big)^2\ra &=& 0 \label{typeI rel1}, \\
\big(1+\la S^*u, v\ra\big)\la u, \big(\frac{z-\alpha}{1-\bar{\alpha}z}\big)^2\ra-\bar{\alpha}^2u(0)&=& 0 \label{typeI rel2},
\end{eqnarray}
and the operator $A$ in \eqref{define operator A} (with $k=1$) is normal. We show that the relations \eqref{typeI rel1}, and \eqref{typeI rel2} are equivalent to the conditions (2), and (3) of the statement and also, the operator $A$ in \eqref{define operator A} (with $k=1$) is a nonzero scalar multiple of the operator $T_1$ as given by $(4)$ of the statement.

Note by \eqref{v su perp kernel}, there exist scalars $c_i, d_i\in\C$, $i=0,1$ with $(c_0, c_1), (d_0, d_1)$ different from $(0,0)$ such that
\begin{eqnarray}
v&=& c_0 k_{\alpha}+c_1 \frac{z-\alpha}{1-\bar{\alpha}z}k_{\alpha}, \label{v typeI} \\
S^*u&=& d_0 k_{\alpha}+d_1 \frac{z-\alpha}{1-\bar{\alpha}z}k_{\alpha} \label{Su typeI}.
\end{eqnarray}
By \eqref{Su typeI}
\begin{equation}\label{u typeI}
u=u(0)+d_0 zk_{\alpha}+d_1 z\Big(\frac{z-\alpha}{1-\bar{\alpha}z}\Big)k_{\alpha},
\end{equation}
and also an easy computation via \eqref{v typeI}, \eqref{Su typeI} yield
\begin{eqnarray}
c_0 &=& (1-|\alpha|^2)v(\alpha), \label{c0 typeI}\\
c_1 &=& (1-|\alpha|^2)B_{\alpha}^*v(\alpha), \label{c1 typeI} \\
d_0 &=& (1-|\alpha|^2)S^*u(\alpha), \label{d0 typeI}\\
d_1 &=& (1-|\alpha|^2)(zB_{\alpha})^*u(\alpha). \label{d1 typeI}
\end{eqnarray}
Note that
\begin{eqnarray}
\Big\la k_{\alpha}, M_z^*\Big(\frac{z-\alpha}{1-\bar{\alpha}z}\Big)^2\Big \ra &=& -\bar{\alpha}. \label{multiplication star k alpha square} \\
\Big\la k_{\alpha}, M_z^*\Big(\frac{z-\alpha}{1-\bar{\alpha}z}\Big)\Big \ra &=& 1. \label{multiplication star k alpha}
\end{eqnarray}
Then
\begin{equation}\label{zv typeI}
\begin{split}
\Big\la zv, \Big(\frac{z-\alpha}{1-\bar{\alpha}z}\Big)^2\Big \ra &= \Big\la c_0zk_{\alpha}+c_1z\frac{z-\alpha}{1-\bar{\alpha}z}k_{\alpha}, \Big(\frac{z-\alpha}{1-\bar{\alpha}z}\Big)^2\Big \ra\quad (\text{by \eqref{v typeI}})\\
&=-\bar{\alpha} c_0 + c_1 \quad (\text{by }\eqref{multiplication star k alpha square}, \text{and }\eqref{multiplication star k alpha})
\end{split}
\end{equation}
Also, it follows by \eqref{v typeI}
\begin{equation}\label{norm v square I}
\|v\|^2=\Big\la c_0k_{\alpha}+c_1\frac{z-\alpha}{1-\bar{\alpha}z}k_{\alpha}, c_0k_{\alpha}+c_1\frac{z-\alpha}{1-\bar{\alpha}z}k_{\alpha}\Big\ra=\frac{1}{1-|\alpha|^2}(|c_0|^2+|c_1|^2),
\end{equation}
and by \eqref{u typeI}, \eqref{multiplication star k alpha square}, and \eqref{multiplication star k alpha}
\begin{equation}\label{u square blaschke}
\Big\la u, \Big(\frac{z-\alpha}{1-\bar{\alpha}z}\Big)^2\Big \ra= \bar{\alpha}^2u(0)-\bar{\alpha}d_0+d_1.
\end{equation}
Now by equations \eqref{zv typeI}, \eqref{norm v square I}, \eqref{u square blaschke}, the relation \eqref{typeI rel1} reduces after a simple computation to
\begin{equation}\label{typeI rel1 reduced}
c_0\Big(\overline{v(\alpha)}(\bar{\alpha}^2u(0)-\bar{\alpha}d_0+d_1)-\bar{\alpha}\Big)
+ c_1\Big(\overline{B_{\alpha}^*v(\alpha)}(\bar{\alpha}^2u(0)-\bar{\alpha}d_0+d_1)+1\Big) =0.
\end{equation}
Next,
\begin{equation}\label{S star u inner prod v}
\begin{split}
\la S^*u, v\ra &=\Big\la d_0k_{\alpha}+ d_1\frac{z-\alpha}{1-\bar{\alpha}z}k_{\alpha},\quad c_0 k_{\alpha}+c_1 \frac{z-\alpha}{1-\bar{\alpha}z}k_{\alpha}\Big\ra \quad(\text{by }\eqref{v typeI}, \eqref{Su typeI})\\
&=\frac{d_0\bar{c_0}}{1-|\alpha|^2}+\frac{d_1\bar{c_1}}{1-|\alpha|^2}\\
&=d_0\overline{v(\alpha)}+d_1\overline{B_{\alpha}^*v(\alpha)}.\quad (\text{by } \eqref{c0 typeI}, \eqref{c1 typeI})
\end{split}
\end{equation}
Now by \eqref{u square blaschke}, \eqref{S star u inner prod v}, equation \eqref{typeI rel2} reduces to
\begin{equation}\label{typeI rel2 reduced}
d_0\Big(\overline{v(\alpha)}(\bar{\alpha}^2u(0)-\bar{\alpha}d_0+d_1)-\bar{\alpha}\Big)
+ d_1\Big(\overline{B_{\alpha}^*v(\alpha)}(\bar{\alpha}^2u(0)-\bar{\alpha}d_0+d_1)+1\Big) =0.
\end{equation}
Note by \eqref{v typeI} and \eqref{Su typeI}, $(c_0, d_0)\neq(0,0)$ and $(c_1, d_1)\neq (0,0)$. Suppose one of $c_0$ or $d_0$ is zero. For the sake of definiteness, let $c_0=0$. Since $v\neq 0$, it follows by \eqref{v typeI}, $c_1\neq 0$. Then by \eqref{typeI rel1 reduced}
\begin{equation}\label{typeI rel1 separate}
\overline{B_{\alpha}^*v(\alpha)}(\bar{\alpha}^2u(0)-\bar{\alpha}d_0+d_1)+1=0,
\end{equation}
which by \eqref{u square blaschke} is equivalent to
\begin{equation}\label{typeI rel1 final}
\overline{B_{\alpha}^*v(\alpha)}\Big\la u, \Big(\frac{z-\alpha}{1-\bar{\alpha}z}\Big)^2\Big \ra+1=0.
\end{equation}
Again, as $d_0$ can not be zero, it follows by \eqref{typeI rel2 reduced} and \eqref{typeI rel1 separate}
\begin{equation}\label{typeI rel2 separate}
\overline{v(\alpha)}(\bar{\alpha}^2u(0)-\bar{\alpha}d_0+d_1)-\bar{\alpha}=0,
\end{equation}
which is again by \eqref{u square blaschke}, equivalent to
\begin{equation}\label{typeI rel2 final}
\overline{v(\alpha)}\Big\la u,\Big(\frac{z-\alpha}{1-\bar{\alpha}z}\Big)^2\Big \ra-\bar{\alpha}=0.
\end{equation}
Similarly, one will obtain same equations \eqref{typeI rel1 separate}---\eqref{typeI rel2 final} if $d_0=0$ or only one of $c_1, d_1$ is zero. Let us now assume that all $c_i, d_i$ are non-zero for $i=1,2$. Then, since $\{v, S^*u\}$ is linearly independent, it follows by \eqref{v typeI} and \eqref{Su typeI}
\begin{equation}\label{coefficint matrix nonzero}
\det \begin{pmatrix}
c_0 & c_1 \\
d_0 & d_1
\end{pmatrix}=(c_0d_1-c_1d_0)\neq0.
\end{equation}
Then, multiplying \eqref{typeI rel1 reduced} by $d_0$ and \eqref{typeI rel2 reduced} by $c_0$, one will have on subtraction
\begin{equation}\label{coeff nonzero same rel1}
(c_1d_0-d_1c_0)\Big(\overline{B_{\alpha}^*v(\alpha)}(\bar{\alpha}^2u(0)-\bar{\alpha}d_0+d_1)+1\Big)=0.
\end{equation}
which is further by \eqref{coefficint matrix nonzero} and \eqref{u square blaschke} becomes
\begin{equation}\label{rel1 deduced coeff nonzero}
\overline{B_{\alpha}^*v(\alpha)}\Big\la u, \Big(\frac{z-\alpha}{1-\bar{\alpha}z}\Big)^2\Big \ra+1=0.
\end{equation}
Again, multiplying \eqref{typeI rel1 reduced} by $d_1$ and \eqref{typeI rel2 reduced} by $c_1$, one will have on subtraction
\begin{equation}\label{coeff nonzero same rel1}
(c_0d_1-d_0c_1)\Big(\overline{v(\alpha)}(\bar{\alpha}^2u(0)-\bar{\alpha}d_0+d_1)-\bar{\alpha}\Big)=0.
\end{equation}
which is further by \eqref{coefficint matrix nonzero} and \eqref{u square blaschke} becomes
\begin{equation}\label{rel2 deduced coeff nonzero}
\overline{v(\alpha)}\Big\la u, \Big(\frac{z-\alpha}{1-\bar{\alpha}z}\Big)^2\Big \ra-\bar{\alpha}=0.
\end{equation}
Hence in all the possible cases of $(c_i,d_i)$, $i=1,2$, one will obtain
\begin{eqnarray}
\overline{v(\alpha)}\Big\la u, \Big(\frac{z-\alpha}{1-\bar{\alpha}z}\Big)^2\Big \ra-\bar{\alpha} &=& 0, \label{rel1 redefined} \\
\overline{B_{\alpha}^*v(\alpha)}\Big\la u, \Big(\frac{z-\alpha}{1-\bar{\alpha}z}\Big)^2\Big \ra+1&=& 0 \label{rel2 redefined}.
\end{eqnarray}
Note by \eqref{rel2 redefined}, $\overline{B_{\alpha}^*v(\alpha)}$ and $\Big\la u, \Big(\frac{z-\alpha}{1-\bar{\alpha}z}\Big)^2\Big\ra$ are both nonzero. Assume $\alpha\neq 0$. Then the equations \eqref{rel1 redefined} and \eqref{rel2 redefined} together imply
\begin{equation}\label{derived rel}
\begin{split}
\Big\la u, \Big(\frac{z-\alpha}{1-\bar{\alpha}z}\Big)^2\Big\ra\big(\overline{v(\alpha)}+\bar{\alpha}\overline{B_{\alpha}^*v(\alpha)}\big)&=0 \\
\iff v(\alpha)+\alpha B_{\alpha}^*v(\alpha)&=0.
\end{split}
\end{equation}
Since $\Big\la u, \Big(\frac{z-\alpha}{1-\bar{\alpha}z}\Big)^2\Big\ra\neq 0$, it follows by \eqref{rel1 redefined}, $v(\alpha)=0$ if and only if $\alpha=0$. Hence the relations \eqref{rel1 redefined}, and \eqref{rel2 redefined} are equivalent to
\begin{eqnarray}
v(\alpha)+\alpha B_{\alpha}^*v(\alpha)&=& 0, \label{rel1 finalize}\\
\overline{B_{\alpha}^*v(\alpha)}\Big\la u, \Big(\frac{z-\alpha}{1-\bar{\alpha}z}\Big)^2\Big \ra+1&=& 0. \label{rel2 finalize}
\end{eqnarray}
Since $T$ is quasinormal (by assumption), by Theorem \ref{rank 2 independent}, the operator $A$ (corresponding to $k=1$) given by \eqref{define operator A} is normal. It follows along the lines of proof of the Theorem \ref{rank 2 independent} that, $A$ is a nonzero scalar multiple of the compression operator $P_1T|_{\text{span}\{v, S^*u\}}$, where $P_1$ is the orthogonal projection of $H^2(\D)$ onto $\text{span}\{v, S^*u\}$. Note by \eqref{v su perp kernel}, $\text{span}\{v, S^*u\}=\text{span}\{k_{\alpha}, \frac{z-\alpha}{1-\bar{\alpha}z}k_{\alpha}\}$. If $P_2$ denotes the orthogonal projection of $H^2(\D)$ onto $\text{span}\{k_{\alpha}, \frac{z-\alpha}{1-\bar{\alpha}z}k_{\alpha}\}$, then the operator $P_2T|_{\text{span}\{k_{\alpha}, \frac{z-\alpha}{1-\bar{\alpha}z}k_{\alpha}\}}=P_1T|_{\text{span}\{v, S^*u\}}$.
We now show that $P_2T|_{\text{span}\{k_{\alpha}, \frac{z-\alpha}{1-\bar{\alpha}z}k_{\alpha}\}}$ is exactly the operator $T_1$ given in the statement. Note that,
$$Tk_{\alpha}=(S+u\otimes v)k_{\alpha}=z k_{\alpha}+\overline{v(\alpha)}u,$$
which by \eqref{u typeI} simplifies to
\begin{equation}\label{T k alpha}
Tk_{\alpha}=\overline{v(\alpha)}u(0)+\big(1+d_0\overline{v(\alpha)}\big)zk_{\alpha}+d_1\overline{v(\alpha)}z\Big(\frac{z-\alpha}{1-\bar{\alpha}z}\Big)k_{\alpha}.
\end{equation}
Next
$$T\Big(\frac{z-\alpha}{1-\bar{\alpha}z}\Big)k_{\alpha}=(S+u\otimes v)\Big(\frac{z-\alpha}{1-\bar{\alpha}z}\Big)k_{\alpha}=z \Big(\frac{z-\alpha}{1-\bar{\alpha}z}\Big)k_{\alpha}+\overline{B_{\alpha}^*v(\alpha)}u,$$

which together with \eqref{u typeI} yield
\begin{equation}\label{T z k alpha}
T\Big(\frac{z-\alpha}{1-\bar{\alpha}z}\Big)k_{\alpha}=\overline{B_{\alpha}^*v(\alpha)}u(0)+d_0\overline{B_{\alpha}^*v(\alpha)}zk_{\alpha}+\big(1+d_1\overline{B_{\alpha}^*v(\alpha)}\big)z\Big(\frac{z-\alpha}{1-\bar{\alpha}z}\Big)k_{\alpha}.
\end{equation}
Note that $1, zk_{\alpha},z\Big(\frac{z-\alpha}{1-\bar{\alpha}z}\Big)k_{\alpha}\perp\Big(\frac{z-\alpha}{1-\bar{\alpha}z}\Big)^2z^n$ for all $n\geq 1$. Hence $$1, zk_{\alpha},z\Big(\frac{z-\alpha}{1-\bar{\alpha}z}\Big)k_{\alpha}\in \text{span }\Big\{k_{\alpha},\Big(\frac{z-\alpha}{1-\bar{\alpha}z}\Big)k_{\alpha}, \Big(\frac{z-\alpha}{1-\bar{\alpha}z}\Big)^2\Big\}.$$

Then a simple computation via \eqref{multiplication star k alpha square}, \eqref{multiplication star k alpha} yield
\begin{eqnarray}
1 &=& (1-|\alpha|^2)k_{\alpha}-\bar{\alpha}(1-|\alpha|^2)\Big(\frac{z-\alpha}{1-\bar{\alpha}z}\Big)k_{\alpha}+\bar{\alpha}^2\Big(\frac{z-\alpha}{1-\bar{\alpha}z}\Big)^2, \label{sub one} \\
zk_{\alpha}&=&\alpha k_{\alpha}+(1-|\alpha|^2)\Big(\frac{z-\alpha}{1-\bar{\alpha}z}\Big)k_{\alpha}-\bar{\alpha}\Big(\frac{z-\alpha}{1-\bar{\alpha}z}\Big)^2, \label{sub z k alpha}\\
z\Big(\frac{z-\alpha}{1-\bar{\alpha}z}\Big)k_{\alpha}&=& \alpha\Big(\frac{z-\alpha}{1-\bar{\alpha}z}\Big)k_{\alpha}+ \Big(\frac{z-\alpha}{1-\bar{\alpha}z}\Big)^2. \label{sub z k alpha square}
\end{eqnarray}
Substituting $1, zk_{\alpha}, z\Big(\frac{z-\alpha}{1-\bar{\alpha}z}\Big)k_{\alpha}$ from \eqref{sub one}, \eqref{sub z k alpha}, \eqref{sub z k alpha square} respectively to \eqref{T k alpha}, a simple computation reveals that
\begin{equation}\label{T k alpha for matrix}
\begin{split}
Tk_{\alpha}&=\Big(\overline{v(\alpha)}u(0)(1-|\alpha|^2)+\alpha(1+d_0\overline{v(\alpha)})\Big)k_{\alpha}+ \\
&\Big(-\bar{\alpha}(1-|\alpha|^2)\overline{v(\alpha)}u(0)+(1-|\alpha|^2)(1+d_0\overline{v(\alpha)})+\alpha d_1\overline{v(\alpha)}\Big)\frac{z-\alpha}{1-\bar{\alpha}z}k_{\alpha}+\\
&\Big(-\bar{\alpha}+\overline{v(\alpha)}(\bar{\alpha}^2 u(0)-\bar{\alpha}d_0+d_1)\Big)\Big(\frac{z-\alpha}{1-\bar{\alpha}z}\Big)^2k_{\alpha}.
\end{split}
\end{equation}
By \eqref{typeI rel2 separate}, the coefficient of $\Big(\frac{z-\alpha}{1-\bar{\alpha}z}\Big)^2k_{\alpha}$ in $Tk_{\alpha}$ above is zero. If $a, b$ denote the coefficients of $k_{\alpha}$ and $\frac{z-\alpha}{1-\bar{\alpha}z}k_{\alpha}$ in \eqref{T k alpha for matrix}, then one can write
\begin{equation}\label{compression T k typeI}
Tk_{\alpha}=a k_{\alpha}+b\frac{z-\alpha}{1-\bar{\alpha}z}k_{\alpha}.
\end{equation}

We now simplify the coefficients $a$ and $b$. Note by \eqref{d0 typeI} the given quantity $R$ is
\begin{equation}\label{define R}
R=d_0-\bar{\alpha}u(0)=(1-|\alpha|^2)S^*u(\alpha)-\bar{\alpha}u(0).
\end{equation}
Then By \eqref{T k alpha for matrix},
$$a=\overline{v(\alpha)}u(0)(1-|\alpha|^2)+\alpha\big(1+d_0\overline{v(\alpha)}\big)=\overline{v(\alpha)}u(0)+\alpha\Big(1+\overline{v(\alpha)}\big(d_0-\bar{\alpha}u(0)\big)\Big)
,$$
and further by \eqref{define R}
\begin{equation}\label{a typeI}
a=\overline{v(\alpha)}u(0)+\alpha\big(1+\overline{v(\alpha)}R\big).
\end{equation}
Next by \eqref{T k alpha for matrix}
\begin{equation}\label{b typeI calculation}
b=-\bar{\alpha}(1-|\alpha|^2)\overline{v(\alpha)}u(0)+(1-|\alpha|^2)(1+d_0\overline{v(\alpha)})+\alpha d_1\overline{v(\alpha)}.
\end{equation}

The equation \eqref{typeI rel2 separate} can also be written as
\begin{equation}\label{sub rel1 separate}
\overline{v(\alpha)}d_1=\overline{\alpha}\Big(1+d_0\overline{v(\alpha)}-\bar{\alpha}\overline{v(\alpha)}u(0)\Big).
\end{equation}
Then substituting $\overline{v(\alpha)}d_1$ from \eqref{sub rel1 separate} to \eqref{b typeI calculation}, a simple computation yield
$$b=1+\overline{v(\alpha)}\big(d_0-\bar{\alpha}u(0)\big),$$
and hence by \eqref{define R}
\begin{equation}\label{b typeI}
b=1+\overline{v(\alpha)}R.
\end{equation}
Similarly, substituting $1, zk_{\alpha}, z\Big(\frac{z-\alpha}{1-\bar{\alpha}z}\Big)k_{\alpha}$ from \eqref{sub one}, \eqref{sub z k alpha}, \eqref{sub z k alpha square} respectively to \eqref{T z k alpha}, a simple computation yield
\begin{equation}\label{T z k for matrix}
\begin{split}
T\Big(\frac{z-\alpha}{1-\bar{\alpha}z}\Big)k_{\alpha} &=\overline{B_{\alpha}^*v(\alpha)}\Big((1-|\alpha|^2)u(0)+\alpha d_0\Big)k_{\alpha}+ \\
&\Big(-\bar{\alpha}(1-|\alpha|^2)u(0)\overline{B_{\alpha}^*v(\alpha)}+(1-|\alpha|^2)d_0\overline{B_{\alpha}^*v(\alpha)}+\alpha\big(1+d_1\overline{B_{\alpha}^*v(\alpha)}\big)\Big)\frac{z-\alpha}{1-\bar{\alpha}z} k_{\alpha}\\
&+\Big(1+ \overline{B_{\alpha}^*v(\alpha)}(\bar{\alpha}^2 u(0)-\bar{\alpha}d_0+d_1)\Big) \Big(\frac{z-\alpha}{1-\bar{\alpha}z}\Big)^2.
\end{split}
\end{equation}
Note by \eqref{typeI rel1 separate}, the coefficient of $\Big(\frac{z-\alpha}{1-\bar{\alpha}z}\Big)^2$ in $T\Big(\frac{z-\alpha}{1-\bar{\alpha}z}\Big)k_{\alpha}$ above is zero. Let $c,d$ denote the coefficients of $k_{\alpha}$, and $\frac{z-\alpha}{1-\bar{\alpha}z}k_{\alpha}$ respectively in the representation of $T\Big(\frac{z-\alpha}{1-\bar{\alpha}z}\Big)k_{\alpha}$ in \eqref{T z k for matrix}. Then it follows by \eqref{typeI rel1 separate} and \eqref{T z k for matrix}
\begin{equation}\label{compression T z k typeI}
T\Big(\frac{z-\alpha}{1-\bar{\alpha}z}\Big)k_{\alpha}=c k_{\alpha}+d \frac{z-\alpha}{1-\bar{\alpha}z}k_{\alpha}.
\end{equation}
Note by \eqref{T z k for matrix},
$$c=\overline{B_{\alpha}^*v(\alpha)}\Big((1-|\alpha|^2)u(0)+\alpha d_0\Big)=\overline{B_{\alpha}^*v(\alpha)}u(0)+\alpha\overline{B_{\alpha}^*v(\alpha)}
\big(d_0-\bar{\alpha}u(0)\big),$$
and hence by \eqref{define R}
\begin{equation}\label{c typeI}
c=\overline{B_{\alpha}^*v(\alpha)}u(0)+\alpha \overline{B_{\alpha}^*v(\alpha)}R.
\end{equation}
Again, by \eqref{typeI rel1 separate}
\begin{equation}\label{sub rel2 separate}
\overline{B_{\alpha}^*v(\alpha)}d_1=\bar{\alpha}\overline{B_{\alpha}^*v(\alpha)}(d_0-\bar{\alpha}u(0))-1.
\end{equation}
By \eqref{T z k for matrix}
\begin{equation}\label{d typeI ist}
d=-\bar{\alpha}(1-|\alpha|^2)u(0)\overline{B_{\alpha}^*v(\alpha)}+(1-|\alpha|^2)d_0\overline{B_{\alpha}^*v(\alpha)}+\alpha\big(1+d_1\overline{B_{\alpha}^*v(\alpha)}\big).
\end{equation}
Now substituting $\overline{B_{\alpha}^*v(\alpha)}d_1$ from \eqref{sub rel2 separate} in \eqref{d typeI ist}, it follows by a simple computation
$$d=\overline{B_{\alpha}^*v(\alpha)}(d_0-\bar{\alpha}u(0)),$$ and hence by \eqref{define R}
\begin{equation}\label{d typeI}
d=\overline{B_{\alpha}^*v(\alpha)}R.
\end{equation}
Therefore, it follows by \eqref{compression T k typeI}, \eqref{compression T z k typeI} \eqref{a typeI}, \eqref{b typeI}, \eqref{c typeI}, and \eqref{d typeI} that the operator $P_2T|_{\text{span}\{k_{\alpha}, \Big(\frac{z-\alpha}{1-\bar{\alpha}z}\Big)k_{\alpha}\}}:\text{span}\{k_{\alpha}, \Big(\frac{z-\alpha}{1-\bar{\alpha}z}\Big)k_{\alpha}\}\longrightarrow \text{span}\{k_{\alpha}, \Big(\frac{z-\alpha}{1-\bar{\alpha}z}\Big)k_{\alpha}\} $ is given by
\begin{equation*}
\begin{split}
P_2T|_{\text{span}\{k_{\alpha}, \Big(\frac{z-\alpha}{1-\bar{\alpha}z}\Big)k_{\alpha}\}}(k_{\alpha})&= \Big(\overline{v(\alpha)}u(0)+\alpha\big(1+\overline{v(\alpha)}R\big)\Big)k_{\alpha}+ \Big(1+\overline{v(\alpha)}R\Big)\frac{z-\alpha}{1-\bar{\alpha}z}k_{\alpha} \\
P_2T|_{\text{span}\{k_{\alpha}, \Big(\frac{z-\alpha}{1-\bar{\alpha}z}\Big)k_{\alpha}\}}\Big(\frac{z-\alpha}{1-\bar{\alpha}z}k_{\alpha}\Big)&=\Big(\overline{B_{\alpha}^*v(\alpha)}u(0)+\alpha \overline{B_{\alpha}^*v(\alpha)}R\Big)k_{\alpha}+\big(\overline{B_{\alpha}^*v(\alpha)}R\big)\frac{z-\alpha}{1-\bar{\alpha}z}k_{\alpha}.
\end{split}
\end{equation*}
Clearly $$P_2T|_{\text{span}\{k_{\alpha}, \Big(\frac{z-\alpha}{1-\bar{\alpha}z}\Big)k_{\alpha}\}}=T_1.$$
Conversely, let $T=S+u\otimes v$, $\{v, S^*u\}$ linearly independent, and the conditions $(1)$---$(4)$ of the statement holds. Since by $(1)$, $v,S^*u\in \text{ span }\{k_{\alpha}, \frac{z-\alpha}{1-\bar{\alpha}z}k_{\alpha}\}$ for some $\alpha\in\D$, it follows that $\{v, S^*u\}$ is $S^*$-invariant and $\{v, S^*u\}^{\perp}=\Big(\frac{z-\alpha}{1-\bar{\alpha}z}\Big)^2H^2(\D)$. Then one can show via condition~(1)
\begin{eqnarray}
v &=& (1-|\alpha|^2)\Big(v(\alpha)k_{\alpha}+B_{\alpha}^*v(\alpha)\frac{z-\alpha}{1-\bar{\alpha}z}k_{\alpha}\Big), \label{v convs} \\
S^*u &=& (1-|\alpha|^2)\Big(S^*u(\alpha) k_{\alpha}+(zB_{\alpha})^*u(\alpha)\frac{z-\alpha}{1-\bar{\alpha}z}k_{\alpha}\Big). \label{S star u convs}
\end{eqnarray}
As we show earlier, $\ker(I-T^*T)=\{v,S^*u\}^{\perp}$, it follows that
$$\ker(I-T^*T)\ominus z\ker(I-T^*T)=\C\Big(\frac{z-\alpha}{1-\bar{\alpha}z}\Big)^2.$$

Now proceeding exactly like the first part, one can show that the conditions $(2), (3)$ of the statement are equivalent to \eqref{rel1},\eqref{rel2} with $g_1(z)=\Big(\frac{z-\alpha}{1-\bar{\alpha}z}\Big)^2$ of the Theorem \ref{rank 2 independent}. Again, since the operator in the statement is normal and can be shown (by the same argument as in the first part) to be the scalar multiple of the operator $A$ in \eqref{define operator A} (with $k=1$), the proof follows by the converse part of the Theorem \ref{rank 2 independent}.
\end{proof}

\subsection{Type II quasinormality}\label{sec: Type 1}

\begin{Theorem}\label{type 2 quasinormal}
Let $T=S+u\otimes v$ be on $H^2(\D)$ with $\{v, S^*u\}$ linearly independent. Then $T$ is quasinormal of type II if and only if there exists distinct $\alpha, \beta\in\D$ such that
\begin{enumerate}
\item $v, S^*u\in\text{ span } \{k_{\alpha}, k_{\beta}\}$,
\item $\overline{v(\alpha)}\Big\la u, \frac{z-\alpha}{1-\bar{\alpha}z}\frac{z-\beta}{1-\bar{\beta}z}\Big\ra-\bar{\beta}=0$,
\item $\overline{v(\beta)}\Big\la u, \frac{z-\alpha}{1-\bar{\alpha}z}\frac{z-\beta}{1-\bar{\beta}z}\Big\ra-\bar{\alpha}=0$,
and
\item If $s_0=\frac{1}{\overline{\alpha-\beta}}\la S^*u, z-\beta\ra$, and $s_1=- \frac{1}{\overline{\alpha-\beta}}\la S^*u, z-\alpha\ra$, then for
$\alpha$ nonzero the operator $T_1: \text{ span } \{k_{\alpha}, k_{\beta}\}\longrightarrow \text{ span } \{k_{\alpha}, k_{\beta}\}$ defined by
\begin{equation}\label{T1 type2}
\begin{split}
T_1(k_{\alpha})&=\frac{1}{\bar{\alpha}}\big(1+\overline{v(\alpha)}s_0\big)k_{\alpha}+ \Big(\overline{v(\alpha)}u(0)-\frac{1}{\bar{\alpha}}\big(1+\overline{v(\alpha)}s_0\big)\Big)k_{\beta} \\
T_1(k_{\beta})&=s_0\frac{\overline{v(\beta)}}{\bar{\alpha}} k_{\alpha}+\Big(\overline{v(\beta)}u(0)-\frac{\overline{v(\beta)}}{\bar{\alpha}}s_0\Big)k_{\beta}
\end{split}
\end{equation}
is normal
or for $\beta$ nonzero, the operator $T_2: \text{ span } \{k_{\alpha}, k_{\beta}\}\longrightarrow \text{ span } \{k_{\alpha}, k_{\beta}\}$ defined by
\begin{equation}\label{T2 type2}
\begin{split}
T_2(k_{\alpha})&=\Big(\overline{v(\alpha)}u(0)-\frac{\overline{v(\alpha)}}{\bar{\beta}}s_1\Big)k_{\alpha}+ s_1\frac{\overline{v(\alpha)}}{\bar{\beta}}k_{\beta}\\
T_2(k_{\beta})&=\Big(\overline{v(\beta)}u(0)-\frac{1}{\bar{\beta}}\big(1+s_1\overline{v(\beta)}\Big)k_{\alpha}+\frac{1}{\bar{\beta}}\Big(1+\overline{v(\beta)}s_1\Big)k_{\beta}
\end{split}
\end{equation}
is normal.
\end{enumerate}
\end{Theorem}

\begin{proof}
Let $T=S+ u \otimes v$, $\{v, S^*u\}$ be linearly independent, and assume $T$ to be quasinormal of type II on $H^2(\D)$. Then $\overline{\text{ran}}(I-T^*T)=\text{span }\{v, S^*u\}$ reduces $T$ and $S^*$-invariant. By the discussion prior to Theorem \ref{type 1 quasinormal}, there exist $\alpha, \beta\in \D$ with $\alpha\neq \beta$ such that
$$\ker(I-T^*T)=\{v, S^*u\}^{\perp}=\Big(\frac{z-\alpha}{1-\bar{\alpha}z}\Big)\Big(\frac{z-\beta}{1-\bar{\beta}z}\Big)H^2(\D).$$

\NI Since $k_{\alpha}, k_{\beta}\perp \big(\frac{z-\alpha}{1-\bar{\alpha}z}\big)\big(\frac{z-\beta}{1-\bar{\beta}z}\big)f$ for all $f\in H^2(\D)$, it follows that $v, S^*u\in \text{ span }\{k_{\alpha}, k_{\beta}\}$, which proves the condition (1). It also follows that there exist scalars $t_i, s_i$, $i=0,1$ with $(t_0,t_1), (s_0,s_1)$ different from $(0,0)$ such that
\begin{eqnarray}
v&=& t_0 k_{\alpha}+t_1 k_{\beta}, \label{v type2}\\
S^*u &=& s_0 k_{\alpha}+s_1 k_{\beta}. \label{S star u type2}
\end{eqnarray}
Note by \eqref{S star u type2}
\begin{equation}\label{u type2}
u=u(0)+s_0zk_{\alpha}+s_1zk_{\beta},
\end{equation}
and a simple computation via \eqref{v type2}, \eqref{S star u type2} implies
\begin{eqnarray}
t_0 &=& \frac{1}{\overline{\alpha-\beta}}\la v, z-\beta\ra, \label{t0 type2} \\
t_1 &=&-\frac{1}{\overline{\alpha-\beta}}\la v, z-\alpha\ra, \label{t1 type2} \\
s_0 &=& \frac{1}{\overline{\alpha-\beta}}\la S^*u, z-\beta\ra, \label{s0 type2}\\
s_1 &=&- \frac{1}{\overline{\alpha-\beta}}\la S^*u, z-\alpha\ra. \label{s1 type2}
\end{eqnarray}
Note that $$\ker(I-T^*T)\ominus z(I-T^*T)=\C \Big(\frac{z-\alpha}{1-\bar{\alpha}z}\Big)\Big(\frac{z-\beta}{1-\bar{\beta}z}\Big).$$

Since $T$ is quasinormal, by Theorem \ref{rank 2 independent}, Equations \eqref{rel1}, \eqref{rel2} hold with $$g_1(z)=\Big(\frac{z-\alpha}{1-\bar{\alpha}z}\Big)\Big(\frac{z-\beta}{1-\bar{\beta}z}\Big),\quad \text{that is}$$
\begin{eqnarray}
\Big\la zv, \big(\frac{z-\alpha}{1-\bar{\alpha}z}\big)\big(\frac{z-\beta}{1-\bar{\beta}z}\big)\Big\ra+ \|v\|^2\Big\la u, \big(\frac{z-\alpha}{1-\bar{\alpha}z}\big)\big(\frac{z-\beta}{1-\bar{\beta}z}\big)\Big\ra &=& 0 \label{type2 rel1},\ \  \\
\big(1+\la S^*u, v\ra\big)\Big\la u, \big(\frac{z-\alpha}{1-\bar{\alpha}z}\big)\big(\frac{z-\beta}{1-\bar{\beta}z}\big)\Big\ra-\overline{\alpha\beta}u(0)&=& 0, \label{type2 rel2}
\end{eqnarray}
\noindent
and the operator $A$ in \eqref{define operator A} (with $k=1$) is normal. We show that \eqref{type2 rel1}, \eqref{type2 rel2} are equivalent to given Equations $(2)$--$(3)$, and the operator $A$ is a nonzero scalar multiple of one of $T_1, T_2$ as defined in the statement.
Notice
\begin{eqnarray}
\Big\la k_{\alpha},\ \ M_z^*\frac{z-\alpha}{1-\bar{\alpha}z}\frac{z-\beta}{1-\bar{\beta}z}\Big\ra &=& -\bar{\beta}, \label{k alpha rel3}\\
\Big\la k_{\beta},\ \ M_z^*\frac{z-\alpha}{1-\bar{\alpha}z}\frac{z-\beta}{1-\bar{\beta}z}\Big\ra &=& -\bar{\alpha}. \label{k beta rel4}
\end{eqnarray}
\begin{equation}\label{ist term rel1 type2}
\Big\la zv, \big(\frac{z-\alpha}{1-\bar{\alpha}z}\big)\big(\frac{z-\beta}{1-\bar{\beta}z}\big)\Big\ra=\Big\la t_0k_{\alpha}+t_1k_{\beta}, \quad M_z^*\frac{z-\alpha}{1-\bar{\alpha}z}\frac{z-\beta}{1-\bar{\beta}z}\Big\ra \quad (\text{by} \eqref{v type2}),
\end{equation}
and further by \eqref{k alpha rel3}, \eqref{k beta rel4}
\begin{equation}\label{ist term final rel1 type2}
\Big\la zv, \big(\frac{z-\alpha}{1-\bar{\alpha}z}\big)\big(\frac{z-\beta}{1-\bar{\beta}z}\big)\Big\ra=-(t_0\bar{\beta}+t_1\bar{\alpha}).
\end{equation}
Again, by \eqref{v type2}
\begin{equation}\label{norm v rel1 type2}
\|v\|^2=\la t_0k_{\alpha}+t_1k_{\beta}, v\ra= t_0\overline{v(\alpha)}+t_1\overline{v(\beta)}.
\end{equation}
Now by \eqref{u type2}
\begin{equation}\label{last term rel 1 type2}
\Big\la u,\ \ \big(\frac{z-\alpha}{1-\bar{\alpha}z}\big)\big(\frac{z-\beta}{1-\bar{\beta}z}\big)\Big\ra=\Big\la u(0)+s_0zk_{\alpha}+s_1zk_{\beta},\ \ \big(\frac{z-\alpha}{1-\bar{\alpha}z}\big)\big(\frac{z-\beta}{1-\bar{\beta}z}\big)\Big\ra,
\end{equation}
which by \eqref{k alpha rel3}, \eqref{k beta rel4}, further reduces to
\begin{equation}\label{last term final rel1 type2}
\Big\la u,\ \ \big(\frac{z-\alpha}{1-\bar{\alpha}z}\big)\big(\frac{z-\beta}{1-\bar{\beta}z}\big)\Big\ra=\overline{\alpha\beta}u(0)-s_0\bar{\beta}-s_1\bar{\alpha}.
\end{equation}
Now by \eqref{ist term final rel1 type2}, \eqref{norm v rel1 type2}, and \eqref{last term final rel1 type2}, the equation \eqref{type2 rel1} simplifies to
\begin{equation}\label{type2 rel1 reduced}
t_0\Big(\overline{v(\alpha)}\big(\overline{\alpha\beta}u(0)-s_0\bar{\beta}-s_1\bar{\alpha}\big)-\bar{\beta}\Big)+t_1\Big(\overline{v(\beta)}\big(\overline{\alpha\beta}u(0)-s_0\bar{\beta}-s_1\bar{\alpha}\big)-\bar{\alpha}\Big)=0.
\end{equation}
Note by \eqref{v type2},
$$\la v, k_{\alpha}\ra= \la t_0k_{\alpha}+t_1k_{\beta}, k_{\alpha}\ra=\frac{t_0}{1-|\alpha|^2}+\frac{t_1}{1-\bar{\beta}\alpha},$$
and hence
\begin{equation}\label{v alpha type2}
v(\alpha)=\frac{t_0}{1-|\alpha|^2}+\frac{t_1}{1-\bar{\beta}\alpha}.
\end{equation}
Similarly it follows by \eqref{v type2}
\begin{equation}\label{v beta type2}
v(\beta)=\la v,k_{\beta}\ra=\frac{t_0}{1-\bar{\alpha}\beta}+\frac{t_1}{1-|\beta|^2}.
\end{equation}
Now by \eqref{v type2}, and \eqref{S star u type2}
\begin{equation}\label{ist term rel2 type2}
\la S^*u, v\ra= \la s_0k_{\alpha}+s_1k_{\beta},\ \ t_0k_{\alpha}+t_1k_{\beta}\ra= s_0\Big(\frac{\bar{t_0}}{1-|\alpha|^2}+\frac{\bar{t_1}}{1-\bar{\alpha}\beta}\Big)+s_1\Big(\frac{\bar{t_0}}{1-\bar{\beta}\alpha}+\frac{\bar{t_1}}{1-|\beta|^2}\Big),
\end{equation}
and further by \eqref{v alpha type2}, and \eqref{v beta type2}
\begin{equation}\label{ist term final rel2 type2}
\la S^*u, v\ra= s_0\overline{v(\alpha)}+s_1\overline{v(\beta)}.
\end{equation}
Substituting \eqref{ist term final rel2 type2}, \eqref{last term final rel1 type2} in \eqref{type2 rel2}, one will obtain after a simplification
\begin{equation}\label{type2 rel2 reduced}
s_0\Big(\overline{v(\alpha)}\big(\overline{\alpha\beta}u(0)-s_0\bar{\beta}-s_1\bar{\alpha}\big)-\bar{\beta}\Big)+s_1\Big(\overline{v(\beta)}\big(\overline{\alpha\beta}u(0)-s_0\bar{\beta}-s_1\bar{\alpha}\big)-\bar{\alpha}\Big)=0.
\end{equation}
Since $v, S^*u$ are linearly independent, it follows by \eqref{v type2}, and \eqref{S star u type2} that both the pair $(t_0, s_0), (t_1, s_1)$ are different from $(0,0)$. There are the possibilities of exactly one element from one or both the pairs $(t_0, s_0)$, $(t_1, s_1)$ is zero or $t_i, s_i$ are nonzero for all $i=1,2$. In the later case, one will have
$$
\det\begin{pmatrix}
t_0 & s_0 \\
t_1 & s_1
\end{pmatrix}
=(t_0s_1-s_0t_1)\neq 0.$$

Then for each of the possible cases, proceeding exactly as in the Theorem \ref{type 1 quasinormal}, one will obtain by \eqref{type2 rel1 reduced}, \eqref{type2 rel2 reduced}
\begin{eqnarray}
\overline{v(\alpha)}\Big(\overline{\alpha\beta}u(0)-s_0\bar{\beta}-s_1\bar{\alpha}\Big)-\bar{\beta}&=& 0, \label{reduced rel1 type2} \\
\overline{v(\beta)}\Big(\overline{\alpha\beta}u(0)-s_0\bar{\beta}-s_1\bar{\alpha}\Big)-\bar{\alpha} &=& 0, \label{reduced rel2 type2}
\end{eqnarray}
which can also be written via \ref{last term final rel1 type2}
\begin{eqnarray}
\overline{v(\alpha)}\Big\la u,\ \ \big(\frac{z-\alpha}{1-\bar{\alpha}z}\big)\big(\frac{z-\beta}{1-\bar{\beta}z}\big)\Big\ra-\bar{\beta}&=& 0, \\
\overline{v(\beta)}\Big\la u,\ \ \big(\frac{z-\alpha}{1-\bar{\alpha}z}\big)\big(\frac{z-\beta}{1-\bar{\beta}z}\big)\Big\ra-\bar{\alpha} &=& 0.
\end{eqnarray}
As we observed, the operator $A$ (equation \eqref{define operator A} with $k=1$) in the Theorem \ref{rank 2 independent} is a nonzero scalar multiple of the operator $P_1T|_{\text{span}\{v, S^*u\}}$ where $P_1$ is the orthogonal projection of $H^2(\D)$ onto $\text{span}\{v, S^*u\}$. Again, since $\text{span}\{k_{\alpha}, k_{\beta}\}=\text{span}\{v, S^*u\}$, it follows that $P_2T|_{\text{span}\{k_{\alpha}, k_{\beta}\}}=P_1T|_{\text{span}\{v, S^*u\}}$, where $P_2$ is the orthogonal projection of $H^2(\D)$ onto $\text{span}\{k_{\alpha}, k_{\beta}\}$. We show that $P_2T|_{\text{span}\{k_{\alpha}, k_{\beta}\}}$ is exactly one of the operators $T_1, T_2$ as defined in the statement. Since $T$ is quasinormal by assumption, the proof of the first part will follow by the forward implication of the Theorem \ref{rank 2 independent}.

Note by \eqref{u type2}
\begin{eqnarray}
Tk_{\alpha}&&=zk_{\alpha}+\overline{v(\alpha)}u=\ \ zk_{\alpha}+\overline{v(\alpha)}\big(u(0)+s_0zk_{\alpha}+s_1zk_{\beta}\big),
 \label{T k alpha type2}\\
Tk_{\beta}&&=zk_{\beta}+\overline{v(\beta)}u=\ \ zk_{\beta}+\overline{v(\beta)}\big(u(0)+s_0zk_{\alpha}+s_1zk_{\beta}
\big).\label{T k beta type2}
\end{eqnarray}
An easy computation shows that
\begin{eqnarray}
\Big\la k_{\alpha},M_z^* \frac{z-\beta}{1-\bar{\beta}z}k_{\alpha}\Big\ra &=& \frac{1-|\beta|^2}{1-\bar{\alpha}\beta}+\frac{\alpha}{1-|\alpha|^2}\frac{\bar{\alpha}-\bar{\beta}}{1-\bar{\alpha}\beta}, \label{1 type2}\\
\Big\la k_{\beta},M_z^* \frac{z-\alpha}{1-\bar{\alpha}z}k_{\beta}\Big\ra &=& \frac{1-|\alpha|^2}{1-\alpha\bar{\beta}}-\frac{\beta}{1-|\beta|^2}\frac{\bar{\alpha}-\bar{\beta}}{1-\alpha\bar{\beta}}, \label{z k alpha type2}\\
\Big\la k_{\alpha}, M_z^* \frac{z-\alpha}{1-\bar{\alpha}z}k_{\beta}\Big\ra &=& 1, \label{k alpha Mz star type2}\\
\Big\la k_{\beta},M_z^* \frac{z-\beta}{1-\bar{\beta}z}k_{\alpha}\Big\ra &=& 1. \label{k beta Mz star type2}
\end{eqnarray}
Since $1, zk_{\alpha}, zk_{\beta}\perp \big(\frac{z-\alpha}{1-\bar{\alpha}z}\big)\big(\frac{z-\beta}{1-\bar{\beta}z}\big)z^n$ for all $n\geq 1$, it follows that
$$ 1, zk_{\alpha}, zk_{\beta}\in\text{span }\{k_{\alpha}, k_{\beta},\big(\frac{z-\alpha}{1-\bar{\alpha}z}\big)\big(\frac{z-\beta}{1-\bar{\beta}z}\big)\},
$$
and it can be shown by \eqref{k alpha rel3}, \eqref{k beta rel4}, and \eqref{1 type2}---\eqref{k beta Mz star type2}
\begin{equation}\label{1 in type2}
1= -\frac{\bar{\beta}(1-\bar{\alpha}\beta)(1-|\alpha|^2)}{\overline{\alpha-\beta}}k_{\alpha}+ \frac{\bar{\alpha}(1-\alpha\bar{\beta})(1-|\beta|^2)}{\overline{\alpha-\beta}}k_{\beta}+\overline{\alpha\beta}\frac{z-\alpha}{1-\bar{\alpha}z}\frac{z-\beta}{1-\bar{\beta}z},
\end{equation}
\begin{equation}\label{z k alpha in type2}
zk_{\alpha}=\Big(\alpha+\frac{(1-|\alpha|^2)(1-|\beta|^2)}{\overline{\alpha-\beta}}\Big)k_{\alpha}-\frac{(1-\bar{\beta}\alpha)(1-|\beta|^2)}{\overline{\alpha-\beta}}k_{\beta}-\bar{\beta}\frac{z-\alpha}{1-\bar{\alpha}z}\frac{z-\beta}{1-\bar{\beta}z},
\end{equation}
\begin{equation}\label{z k beta in type2}
zk_{\beta} = \frac{(1-\bar{\alpha}\beta)(1-|\alpha|^2)}{\overline{\alpha-\beta}}k_{\alpha}+\Big(\beta-\frac{(1-|\alpha|^2)(1-|\beta|^2)}{\overline{\alpha-\beta}}\Big)k_{\beta}-\bar{\alpha}\frac{z-\alpha}{1-\bar{\alpha}z}\frac{z-\beta}{1-\bar{\beta}z}.
\end{equation}
Substituting  $1, zk_{\alpha}, zk_{\beta}$ from \eqref{1 in type2}---\eqref{z k beta in type2} in \eqref{T k alpha type2}, a simplification implies
\begin{equation}\label{T k alpha type2 simpli}
\begin{split}
Tk_{\alpha} &= \Big[(1+\overline{v(\alpha)}s_0)\Big(\alpha+\frac{(1-|\alpha|^2)(1-|\beta|^2)}{\overline{\alpha-\beta}}\Big)+\overline{v(\alpha)}(s_1-\bar{\beta}u(0))\frac{(1-\bar{\alpha}\beta)(1-|\alpha|^2)}{\overline{\alpha-\beta}}\Big]k_{\alpha}+ \\
&\Big[\frac{(1-\bar{\beta}\alpha)(1-|\beta|^2)}{\overline{\alpha-\beta}}\Big(-1+\overline{v(\alpha)}(\bar{\alpha}u(0)-s_0)\Big)+\overline{v(\alpha)}s_1\big(\beta-\frac{(1-|\alpha|^2)(1-|\beta|^2)}{\overline{\alpha-\beta}}\big)\Big]k_{\beta} \\
&\Big[\overline{v(\alpha)}\Big(\overline{\alpha\beta}u(0)-s_0\bar{\beta}-s_1\bar{\alpha}\Big)-\bar{\beta}\Big]\frac{z-\alpha}{1-\bar{\alpha}z}\frac{z-\beta}{1-\bar{\beta}z}.
\end{split}
\end{equation}
By \eqref{reduced rel1 type2}, the coefficient of $\frac{z-\alpha}{1-\bar{\alpha}z}\frac{z-\beta}{1-\bar{\beta}z}$  in $Tk_{\alpha}$ above is zero. Suppose $c_{\alpha}, d_{\alpha}$ denote the coefficients of $k_{\alpha}, k_{\beta}$ respectively in T$k_{\alpha}$ above. Then the equation \eqref{T k alpha type2 simpli} can be written as
\begin{equation}\label{T k alpha type 2 short}
Tk_{\alpha}=c_{\alpha}k_{\alpha}+d_{\alpha}k_{\beta}.
\end{equation}
Note that, \eqref{reduced rel1 type2} can also be written as
\begin{equation}\label{reduced rel1 type2 use}
-\bar{\alpha}\overline{v(\alpha)}\big(s_1-\bar{\beta}u(0)\big)=\bar{\beta}\big(1+\overline{v(\alpha)}s_0\big).
\end{equation}
Note that, one of $\alpha$ and $\beta$ must be nonzero (as $\alpha\neq\beta$ by assumption). Let $\alpha\neq 0$. Then by \eqref{T k alpha type2 simpli}, \eqref{reduced rel1 type2 use}, $c_\alpha$ in \eqref{T k alpha type 2 short} reduces to
\begin{equation}\label{c alpha}
\begin{split}
c_{\alpha}&=(1+\overline{v(\alpha)}s_0)\Big(\alpha+\frac{(1-|\alpha|^2)(1-|\beta|^2)}{\overline{\alpha-\beta}}\Big)-\frac{\bar{\beta}}{\bar{\alpha}}(1+\overline{v(\alpha)}s_0)\frac{(1-\bar{\alpha}\beta)(1-|\alpha|^2)}{\overline{\alpha-\beta}} \\
&=(1+\overline{v(\alpha)}s_0)\Big[\alpha+\frac{(1-|\alpha|^2)(1-|\beta|^2)}{\overline{\alpha-\beta}}-\frac{\bar{\beta}}{\bar{\alpha}}\frac{(1-\bar{\alpha}\beta)(1-|\alpha|^2)}{\overline{\alpha-\beta}}\Big],
\end{split}
\end{equation}
and a further simplification yields
\begin{equation}\label{c alpha nonzero}
c_{\alpha}=\frac{1}{\bar{\alpha}}\big(1+\overline{v(\alpha)}s_0\big).
\end{equation}
If $\beta\neq 0$, then by \eqref{T k alpha type2 simpli}, and \eqref{reduced rel1 type2 use}, $c_{\alpha}$ in \eqref{T k alpha type 2 short} becomes
\begin{equation}\label{c alpha for beta}
\begin{split}
c_{\alpha}&=-\frac{\bar{\alpha}}{\bar{\beta}}\overline{v(\alpha)}(s_1-\bar{\beta}u(0))\Big(\alpha+\frac{(1-|\alpha|^2)(1-|\beta|^2)}{\overline{\alpha-\beta}}\Big)+\overline{v(\alpha)}(s_1-\bar{\beta}u(0))\frac{(1-\bar{\alpha}\beta)(1-|\alpha|^2)}{\overline{\alpha-\beta}} \\
&=\overline{v(\alpha)}(s_1-\bar{\beta}u(0))\Big[-\frac{\bar{\alpha}}{\bar{\beta}}\big(\alpha+\frac{(1-|\alpha|^2)(1-|\beta|^2)}{\overline{\alpha-\beta}}\big)+\frac{(1-\bar{\alpha}\beta)(1-|\alpha|^2)}{\overline{\alpha-\beta}}\Big],
\end{split}
\end{equation}
and further by a simple computation one will have
\begin{equation}\label{c alpha last}
c_{\alpha}=-\frac{\overline{v(\alpha)}}{\bar{\beta}}\big(s_1-\bar{\beta}u(0)\big).
\end{equation}
We now simplify $d_{\alpha}$ in \eqref{T k alpha type 2 short}. Note by \eqref{reduced rel1 type2 use}
\begin{equation}\label{reduced rel1 type 2 2nd use}
\bar{\alpha}\overline{v(\alpha)}s_1=\bar{\beta}\Big[-1+\overline{v(\alpha)}\big(\bar{\alpha}u(0)-s_0\big)\Big].
\end{equation}
If $\alpha\neq 0$, the one can write via \eqref{T k alpha type2 simpli}, and \eqref{reduced rel1 type 2 2nd use}
\begin{equation}\label{d alpha}
d_{\alpha}=\big(-1+\overline{v(\alpha)}(\bar{\alpha}u(0)-s_0)\big)\Big[\frac{(1-\bar{\beta}\alpha)(1-|\beta|^2)}{\overline{\alpha-\beta}}+\frac{\bar{\beta}}{\bar{\alpha}}(\beta-\frac{(1-|\alpha|^2)(1-|\beta|^2)}{\overline{\alpha-\beta}}\big)\Big],
\end{equation}
and further by a simplification
\begin{equation}\label{d alpha nonzero}
d_{\alpha}=\overline{v(\alpha)}u(0)-\frac{1}{\bar{\alpha}}\big(1+\overline{v(\alpha)}s_0\big).
\end{equation}
Again, if $\beta\neq 0$, then it follows by \eqref{T k alpha type2 simpli}, and \eqref{reduced rel1 type 2 2nd use}
\begin{equation}\label{d alpha beta nonzero}
d_{\alpha}=\overline{v(\alpha)}s_1\Big[\frac{\bar{\alpha}}{\bar{\beta}}\frac{(1-\bar{\beta}\alpha)(1-|\beta|^2)}{\overline{\alpha-\beta}}+\beta-\frac{(1-|\alpha|^2)(1-|\beta|^2)}{\overline{\alpha-\beta}}\Big],
\end{equation}
and a further simplification yields
\begin{equation}\label{d alpha final}
d_{\alpha}=\frac{\overline{v(\alpha)}}{\bar{\beta}}s_1.
\end{equation}
Hence it follows by \eqref{T k alpha type 2 short}, \eqref{c alpha nonzero}, \eqref{d alpha nonzero}, \eqref{c alpha last}, \eqref{d alpha final}
\begin{equation}\label{T k alpha final}
Tk_{\alpha}=\begin{cases}
\frac{1}{\bar{\alpha}}\big(1+\overline{v(\alpha)}s_0\big)k_{\alpha}+\Big(\overline{v(\alpha)}u(0)-\frac{1}{\bar{\alpha}}\big(1+\overline{v(\alpha)}s_0\big)\Big)k_{\beta}, & \mbox{if } \alpha\neq 0 \\
-\frac{\overline{v(\alpha)}}{\bar{\beta}}\big(s_1-\bar{\beta}u(0)\big)k_{\alpha}+\frac{\overline{v(\alpha)}}{\bar{\beta}}s_1 k_{\beta}, & \mbox{if } \beta\neq 0.
\end{cases}
\end{equation}
Next, substituting  $1, zk_{\alpha}, zk_{\beta}$ from \eqref{1 in type2}---\eqref{z k beta in type2} in \eqref{T k beta type2} a simplification implies
\begin{equation}\label{T k beta type2 simpli}
\begin{split}
Tk_{\beta} &= \Big[\frac{(1-\bar{\alpha}\beta)(1-|\alpha|^2)}{\overline{\alpha-\beta}}\Big(1+\overline{v(\beta)}\big(s_1-\bar{\beta}u(0)\big)\Big)+\overline{v(\beta)}s_0\Big(\alpha+\frac{(1-|\alpha|^2)(1-|\beta|^2)}{\overline{\alpha-\beta}}\Big)\Big]k_{\alpha}\\
&\Big[\frac{(1-\bar{\beta}\alpha)(1-|\beta|^2)}{\overline{\alpha-\beta}}\Big(\overline{v(\beta)}(\bar{\alpha}u(0)-s_0)\Big)+(1+\overline{v(\beta)}s_1)\big(\beta-\frac{(1-|\alpha|^2)(1-|\beta|^2)}{\overline{\alpha-\beta}}\big)\Big]k_{\beta} \\
&\Big[\overline{v(\beta)}\Big(\overline{\alpha\beta}u(0)-s_0\bar{\beta}-s_1\bar{\alpha}\Big)-\bar{\alpha}\Big]\frac{z-\alpha}{1-\bar{\alpha}z}\frac{z-\beta}{1-\bar{\beta}z}.
\end{split}
\end{equation}
By \eqref{reduced rel2 type2}, the coefficient of $\frac{z-\alpha}{1-\bar{\alpha}z}\frac{z-\beta}{1-\bar{\beta}z}$ in $Tk_{\beta}$ above vanishes. Let $c_{\beta}, d_{\beta}$ denote the coefficients of $k_{\alpha}$ and $k_{\beta}$ respectively in $Tk_{\beta}$ above. Then \eqref{T k beta type2 simpli} becomes
\begin{equation}\label{T k beta short}
Tk_{\beta}=c_{\beta}k_{\alpha}+d_{\beta}k_{\beta}.
\end{equation}
We now simplify $c_{\beta}, d_{\beta}$. Note that \eqref{reduced rel2 type2} can be written as
\begin{equation}\label{rel2 use}
-\bar{\alpha}\Big[1+\overline{v(\beta)}(s_1-\bar{\beta}u(0))\Big]=\bar{\beta}\overline{v(\beta)}s_0.
\end{equation}
If $\alpha\neq 0$, one will have by \eqref{rel2 use} and \eqref{T k beta type2 simpli}
\begin{equation}\label{c beta}
c_{\beta} =-\frac{\bar{\beta}}{\bar{\alpha}}\overline{v(\beta)}s_0 \frac{(1-\bar{\alpha}\beta)(1-|\alpha|^2)}{\overline{\alpha-\beta}}+s_0\overline{v(\beta)}\Big(\alpha+\frac{(1-|\alpha|^2)(1-|\beta|^2)}{\overline{\alpha-\beta}}\Big),\\
\end{equation}
which further reduced to (by a simple computation)
\begin{equation}\label{c beta final}
c_{\beta}=\frac{\overline{v(\beta)}}{\bar{\alpha}}s_0.
\end{equation}
Again, \eqref{reduced rel2 type2} can also be written as
\begin{equation}\label{rel2 2nd use}
\bar{\alpha}\Big(1+\overline{v(\beta)}s_1\Big)=\bar{\beta}\overline{v(\beta)}\Big(\bar{\alpha}u(0)-s_0\Big).
\end{equation}
Hence for $\alpha\neq 0$, it follows via \eqref{rel2 2nd use}, and \eqref{T k beta type2 simpli}
\begin{equation}\label{d beta}
d_{\beta}=\frac{(1-\bar{\beta}\alpha)(1-|\beta|^2)}{\overline{\alpha-\beta}}\overline{v(\beta)}\Big(\bar{\alpha}u(0)-s_0\Big)+\frac{\bar{\beta}}{\bar{\alpha}}\overline{v(\beta)}\Big(\bar{\alpha}u(0)-s_0\Big)\Big(\beta-\frac{(1-|\alpha|^2)(1-|\beta|^2)}{\overline{\alpha-\beta}}\Big),
\end{equation}
and again simplifies to
\begin{equation}\label{d beta final}
d_{\beta}=\frac{\overline{v(\beta)}}{\bar{\alpha}}\big(\bar{\alpha}u(0)-s_0\big).
\end{equation}
Next, if $\beta\neq 0$, it follows by \eqref{T k beta type2 simpli}, and \eqref{rel2 use}
\begin{equation}\label{c beta 2nd use}
c_{\beta}=\frac{(1-\bar{\alpha}\beta)(1-|\alpha|^2)}{\overline{\alpha-\beta}}\Big(1+\overline{v(\beta)}\big(s_1-\bar{\beta}u(0)\big)\Big)-\frac{\bar{\alpha}}{\bar{\beta}}\Big(1+\overline{v(\beta)}\big(s_1-\bar{\beta}u(0)\big)\Big)\Big(\alpha+\frac{(1-|\alpha|^2)(1-|\beta|^2)}{\overline{\alpha-\beta}}\Big),
\end{equation}
and a further simplification yield
\begin{equation}\label{c beta last}
c_{\beta}=\overline{v(\beta)}u(0)-\frac{1}{\bar{\beta}}\Big(1+s_1\overline{v(\beta)}\Big).
\end{equation}
Now for $\beta\neq 0$, it follows by \eqref{T k beta type2 simpli}, and \eqref{rel2 2nd use}
\begin{equation}\label{d beta 2nd use}
d_{\beta}=\frac{\bar{\alpha}}{\bar{\beta}}\frac{(1-\bar{\beta}\alpha)(1-|\beta|^2)}{\overline{\alpha-\beta}}(1+\overline{v(\beta)}s_1)+(1+\overline{v(\beta)}s_1)\big(\beta-\frac{(1-|\alpha|^2)(1-|\beta|^2)}{\overline{\alpha-\beta}}\big),
\end{equation}
and further simplifies to
\begin{equation}\label{d beta last}
d_{\beta}=\frac{1}{\bar{\beta}}\Big(1+\overline{v(\beta)}s_1\Big).
\end{equation}
Therefore it follows by \eqref{T k beta short}, \eqref{c beta final}, \eqref{d beta final}, \eqref{c beta last}, \eqref{d beta last}
\begin{equation}\label{T k beta final}
Tk_{\beta}=\begin{cases}
\frac{\overline{v(\beta)}}{\bar{\alpha}}s_0 k_{\alpha}+\frac{\overline{v(\beta)}}{\bar{\alpha}}\big(\bar{\alpha}u(0)-s_0\big)k_{\beta}, & \mbox{if } \alpha\neq 0 \\
\Big(\overline{v(\beta)}u(0)-\frac{1}{\bar{\beta}}\big(1+s_1\overline{v(\beta)}\big)\Big) k_{\alpha}+\frac{1}{\bar{\beta}}\Big(1+\overline{v(\beta)}s_1\Big)k_{\beta}, & \mbox{if } \beta\neq 0.
\end{cases}
\end{equation}
Hence it follows by \eqref{T k alpha final}, and \eqref{T k beta final} that the operator
$$P_2T|_{\text{span}\{k_{\alpha}, k_{\beta}\}}: \text{span}\{k_{\alpha}, k_{\beta}\} \longrightarrow \text{span}\{k_{\alpha}, k_{\beta}\}$$ given by
\begin{equation*}
\begin{split}
P_2T|_{\text{span}\{k_{\alpha}, k_{\beta}\}}(k_{\alpha})&=\frac{1}{\bar{\alpha}}\big(1+\overline{v(\alpha)}s_0\big)k_{\alpha}+ \Big(\overline{v(\alpha)}u(0)-\frac{1}{\bar{\alpha}}\big(1+\overline{v(\alpha)}s_0\big)\Big)k_{\beta} \\
P_2T|_{\text{span}\{k_{\alpha}, k_{\beta}\}}(k_{\beta})&=s_0\frac{\overline{v(\beta)}}{\bar{\alpha}} k_{\alpha}+\Big(\overline{v(\beta)}u(0)-\frac{\overline{v(\beta)}}{\bar{\alpha}}s_0\Big)k_{\beta}
\end{split}
\end{equation*}
is normal if $\alpha\neq 0$ or
\begin{equation*}
\begin{split}
P_2T|_{\text{span}\{k_{\alpha}, k_{\beta}\}}(k_{\alpha})&=\Big(\overline{v(\alpha)}u(0)-\frac{\overline{v(\alpha)}}{\bar{\beta}}s_1\Big)k_{\alpha}+ s_1\frac{\overline{v(\alpha)}}{\bar{\beta}} k_{\beta}\\
P_2T|_{\text{span}\{k_{\alpha}, k_{\beta}\}}(k_{\beta})&=\Big(\overline{v(\beta)}u(0)-\frac{1}{\bar{\beta}}\big(1+s_1\overline{v(\beta)}\Big)k_{\alpha}+\frac{1}{\bar{\beta}}\Big(1+\overline{v(\beta)}s_1\Big)k_{\beta}
\end{split}
\end{equation*}
is normal if $\beta$ is nonzero.

Conversely let $T=S+u\otimes v$ on $H^2(\D)$ with $\{v, S^*u\}$ linearly independent, and let there be distinct $\alpha,\beta\in\D$ satisfying conditions $(1)$---$(4)$ of the statement. Since by $(1)$, $v, S^*u\in\text{span }\{k_{\alpha}, k_{\beta}\}$, it follows that $\{v, S^*u\}$ is $S^*$-invariant, and also $\{v, S^*u\}^{\perp}=\Big(\frac{z-\alpha}{1-\bar{\alpha}z}\Big)\Big(\frac{z-\beta}{1-\bar{\beta}z}\Big)H^2(\D)$. We have shown earlier, $\ker(I-T^*T)=\{v,S^*u\}^{\perp}$ (case 2, Lemma \ref{Ker-rank1}), and hence
$$\ker(I-T^*T)\ominus z\ker(I-T^*T)=\C\frac{z-\alpha}{1-\bar{\alpha}z}\frac{z-\beta}{1-\bar{\beta}z}.$$
By $(1)$, it also follows that there exist scalars $t_i, s_i\in\C$ such that
$$v= t_0 k_{\alpha}+t_1 k_{\beta},\quad \text{and} \quad S^*u = s_0 k_{\alpha}+s_1 k_{\beta},$$ where the scalars $t_i, s_i$ are given by the same equations \eqref{t0 type2}---\eqref{s1 type2}.

Now proceeding exactly as in the first part, one can show, conditions $(2)$--$(3)$ of the statement are equivalent to Equations \eqref{rel1} and \eqref{rel2} with $k=1$, and $g_1(z)=\frac{z-\alpha}{1-\bar{\alpha}z}\frac{z-\beta}{1-\bar{\beta}z}$ of Theorem \eqref{rank 2 independent}. Similarly, as in the first part, the operators $T_1, T_2$ (for $\alpha\neq 0$ and $\beta\neq 0$) given by $(4)$ in the statement are nonzero scalar multiples of the operator $A$ in \eqref{define operator A} with $k=1$. Since they are normal by the same condition $(4)$, the proof follows by the converse of Theorem \ref{rank 2 independent}.
\end{proof}
\begin{Remark}
Theorems \ref{type 1 quasinormal}, and \ref{type 2 quasinormal} also provide counterexamples to Proposition 2.5 in \cite{Ko-Lee}.
\end{Remark}
\textbf{Acknowledgement:}
The author is thankful to Prof. E. K. Narayanan for many helpful discussions, and suggestions throughout the work. The research of the author is supported by the NBHM-Postdoctoral fellowship at the Indian Institute of Science, India.
\medskip

\textbf{AI Declaration:} The Statements and the proofs of the work are not generated by AI.
\medskip


\bibliographystyle{amsplain}

\begin{thebibliography}{99}

\bibitem{A. Brown}
A. Brown, {\em On a class of operators}, Proc. Amer. Math. Soc., 4 (1953), 723-728.

\bibitem{A. Chattopadhyay A. Skripka}
A. Chattopadhyay, A. Skripka, {\em Trace formulas for relative Schatten class perturbations} J. Funct. Anal., 274 (2018) 3377–3410.

\bibitem{Das}
S. Das, {\em Completely non-unitary contractions and analyticity}, J. Oper. Theory, 94 (2025), 65--91.

\bibitem {Das 1}
S. Das, {\em Invariant subspaces of compression of the Hardy shift on some parametric spaces}, Studia Math., 287 (2026), 21-55.

\bibitem{Das-Sarkar 1}
S. Das, J. Sarkar, {\em Left-Invertibility of Rank-One Perturbations}, Complex Anal. Oper. Theory, 16:109 (2022), 361–382.

\bibitem{Das-Sarkar}
S. Das, J. Sarkar, {\em Invariant subspaces of analytic perturbations}, St. Petersburg Math. J, 35 (2024), 677--695.


\bibitem{Paulsen Raghupati}
K. Davidson, V. Paulsen, M. Raghupathi, D. Singh, {\em A constrained Nevanlinna-Pick interpolation problem}, Indiana Univ. Math. J., 58 (2009), 709-732.

\bibitem{Douglas}
R. G. Douglas, {\em On Majorization, Factorization, and Range Inclusion of Operators on Hilbert Space}, Proc. Amer. Math. Soc., 17 (1966), 413--415.

\bibitem{Foias-Jung-Ko-Percy 1}
C. Foias, I. B. Jung, E. Ko, C. Percy, {\em On rank-one perturbations of normal operators}, J. Funct. Anal. 253 (2007), 628–646.

\bibitem{Foias-Jung-Ko-Percy 2}
C. Foias, I. B. Jung, E. Ko, C. Percy, {\em On rank-one perturbations of normal operators II}, Indiana Univ. Math. J. 57 (2008), 2745–2760.

\bibitem{Halmos}
P. Halmos, {\em A Hilbert space problem book}. Second edition. Graduate Texts in Mathematics, 19. Encyclopedia of Mathematics and its Applications.
17. Springer-Verlag, New York-Berlin, 1982.

\bibitem{E. Ionascu}
E. J. Ionascu, {\em Rank-one perturbations of diagonal operators}, Integr. Equat. Oper. Th., 39 (2001),421–440.

\bibitem{Jung-Lee}
I. B. Jung, E. Y. Lee, {\em Rank one perturbations of Normal operators and hyponormality}, Oper. matrices, 8 (2014), 691–698.

\bibitem{Ko-Lee}
E. Ko, J. E. Lee, {\em On rank one perturbations of the unilateral shift}, Commun. Kor. Math. Soc., 26 (2011), 79–88.

\bibitem{Martin-Putinar}
M. Martin and M. Putinar, {\em Lectures on hyponormal operators, Operator Theory: Advances and Applications}, Vol. 39, Birkhäuser Verlag, Basel, 1989.

\bibitem{NF}
B. Sz.-Nagy, C. Foias, {\em Harmonic analysis of operators on Hilbert space}, North Holland, Amsterdam, 1970.

\bibitem{Nakamura 1}
Y. Nakamura, {\em One-dimensional perturbations of isometries}, Integr. Equat. Oper. Th., 9 (1986), 286--294.

\bibitem{Nakamura}
Y. Nakamura, {\em One-dimensional perturbations of the shift}, Integr. Equat. Oper. Th., 17 (1993), 373--403.

\bibitem{Serban-Turcu}
I. Serban, F. Turcu, {\em Compact perturbations of isometries}, Proc. Amer. Math. Soc., 135, 4, (2007) 1175--1180.
\end{thebibliography}

\end{document}